\documentclass[12pt]{article}
\usepackage{latexsym}
\usepackage{amsmath}
\usepackage{amssymb}
\usepackage[dvips]{graphicx}
\usepackage{latexsym}

\newcommand{\C}{\mathbb{C}}
\newcommand{\CP}{\mathbb{CP}}

\newcommand{\R}{\mathbb{R}}
\newcommand{\Z}{\mathbb{Z}}

\newcommand{\hook}{{\setlength{\unitlength}{11pt}   
                   \begin{picture}(.833,.8)
                   \put(.15,.08){\line(1,0){.35}}
                   \put(.5,.08){\line(0,1){.5}}
                   \end{picture}}}

\def\p{\partial}
\def\ll{\lambda}

\def\be{\begin{equation}}
\def\ee{\end{equation}}
\def\a{\alpha}
\def\lra{\longrightarrow}
\def\ra{\rightarrow}
\def\t{\tilde}
\def\h{\hat}
\def\b{\bar}
\def\om{\omega}
\newcommand{\koniec}{\begin{flushright}  $\Box $ \end{flushright}}
\newtheorem{theo}{Theorem}[section] 
\newtheorem{prop}[theo]{Proposition}  
\newtheorem{lemma}[theo]{Lemma}

\begin{document}
\title{\vskip -70pt
\begin{flushright}
{\normalsize DAMTP-2008-84} \\
\end{flushright}
\vskip 80pt
{\bf Strominger--Yau--Zaslow geometry, Affine Spheres 
and Painlev\'e III.\vskip 20pt}}
\author{Maciej Dunajski\thanks{email M.Dunajski@damtp.cam.ac.uk}\\[15pt]
and\\[15pt]
Prim Plansangkate\thanks{email P.Plansangkate@damtp.cam.ac.uk}
\\[10pt]
{\sl Department of Applied Mathematics and Theoretical Physics} \\[5pt]
{\sl University of Cambridge} \\[5pt]
{\sl Wilberforce Road, Cambridge CB3 0WA, UK}
}
\date{}
\maketitle
\begin{abstract}
{We give a gauge invariant characterisation of the elliptic 
affine sphere equation and the closely related Tzitz\'eica equation
as reductions of real forms of $SL(3, \C)$ anti--self--dual Yang--Mills  equations  by two translations, or equivalently as a special case
of the Hitchin equation.

We use the Loftin--Yau--Zaslow construction to 
give an explicit expression for a six--real dimensional
semi--flat Calabi--Yau metric in terms of a solution
to the affine-sphere equation and show how a subclass
of such metrics arises from 3rd Painlev\'e transcendents.
}
\end{abstract}
\newpage

\section{Introduction}
\setcounter{equation}{0}
Let $X$ be a six real dimensional Calabi--Yau (CY) manifold - a complex 
K\"ahler 
three-fold with covariantly constant holomorphic three-form $\Omega$.  
Any such manifold admits a Ricci flat K\"ahler metric
with holonomy contained in $SU(3)$.

We shall consider a subclass of CY manifolds which are fibred
over a real three dimensional manifold $B$, and the fibres are
special Lagrangian tori $T^3$. This  means that there exists a projection
\[
\pi:X\longrightarrow B
\]
such that the restrictions  of the K\"ahler form $\omega$ and the real part of
the holomorphic three-form $\mbox{Re}(\Omega)$ vanish on
 any fibre $\pi^{-1}(p)\cong T^3$ over a point $p \in B.$ 

The corresponding CY metric is called semi--flat if it is flat
along the fibres. Consider
the K\"ahler form $\omega=i\p\overline{\p}\phi$, where $\phi$
is the K\"ahler potential.  A natural class of semi--flat CY manifolds
are the $T^3$ invariant manifolds. In this case the potential $\phi$
can be chosen not to depend on the coordinates of the fibres
of $\pi.$ The Ricci--flat condition 
$\det \left(\frac{\p^2 \phi}{\p z^j \p \b z^k}\right)=1$ then reduces to the real Monge--Amp\'ere equation
\be
\label{hessian}
\mbox{det}\Big(\frac{\p^2 \phi}{\p x^j \p x^k}\Big)=1,
\ee
where $x^j, j=1,2,3,$ are local coordinates on $B$.
The work of Cheng and Yau \cite{Cheng_Yau} shows that semi--flat
CY metrics on compact complex three-fold are flat, 
so  in what follows we  allow CY manifolds to be non--compact, 
and some fibres of $\pi$ to be singular.

The conjecture of Strominger, Yau and Zaslow  (SYZ) \cite{SYZ}
states that
near the large complex structure limit both $X$ and its mirror
should be the fibrations over
the moduli space of special Lagrangian tori. 
More precisely, SYZ consider the moduli space of special Lagrangian 
submanifolds admitting a unitary
flat connection. They write down a metric on $X$ 
and compute the metric on the moduli space. In the tree level
contribution this metric is derived from the Born--Infeld
action for the brane, assuming that the moduli parameters
slowly vary in time and expanding the action up to second order in time derivatives.
The metric on the moduli space $Y$ arises
from the kinetic term in the Born--Infeld action.
This method is based on Manton's moduli space approximation \cite{manton} and 
was originally used by SYZ. The metric resulting on $Y$ admits
the $T^3$ action even if the original metric on $X$ does not. 
The full agreement between $Y$ and the mirror of $X$ is therefore
expected when instanton contribution from minimal area holomorphic discs
whose boundaries wrap the tori are taken into account.
These corrections are suppressed in the large complex structure limit.

One approach to a proof of the 
Strominger Yau Zaslow conjecture \cite{SYZ}
would be to describe Ricci-flat metrics on Calabi-Yau manifolds near 
large complex structure limits. 
It is 
expected that in the large complex structure limit the base of the
fibration $\pi:X\longrightarrow B$ admits an
affine structure and a special metric of Hessian form.
To test this conjecture Loftin, Yau and Zaslow (LYZ) \cite{LYZ}
aimed to prove the existence of the metric of  
Hessian form\footnote{It follows from the work of Hitchin \cite{Hitchin} that the natural
Weil-Petersson metric on the space of special Lagrangian 
submanifolds has this form.
More precisely, it is shown in \cite{Hitchin} that the 
K\"ahler potentials of $X$ and its mirror $Y$ both satisfy
the Monge-Amp\'ere equation (\ref{hessian}) and are related by a Legendre
transform on the  base. The fibres of the special Lagrangian
fibration of $Y$ are dual (by a Fourier transform)  tori to the fibres 
of $\pi:X\longrightarrow B$.}
\be
\label{Hessian_met}
g_B=\frac{\p^2 \phi}{\p x^j\p x^k} dx^j\otimes dx^k,
\ee
where
$\phi$ is homogeneous of degree 2 in $x^j$ and satisfies (\ref{hessian}). Given such a Hessian metric on $B$, the semi--flat
Calabi--Yau metric $g$ on $TB$ and the corresponding K\"ahler form 
are given by
\be
\label{g_omega}
g=\phi_{jk}(d x^j\otimes d x^k+d y^j\otimes d y^k),\qquad
\omega=\frac{i}{2}\phi_{jk} d z^j \wedge d \overline{z}^k,
\ee
where $y^j$ are coordinates on the fibres of $TB$ and
$z^j=x^j+iy^j$.

LYZ constructed a candidate for such 
metric as a cone over the elliptic affine sphere metric 
with three singular points.  One consequence of Mirror Conjecture is
that the base metric $g_B$ should have singularities in codimension
two, and LYZ were interested in a local metric model near the trivalent
vertex of a Y-shaped singularity.
The monodromy of 
the resulting affine structure has not been calculated, so it is not yet 
clear that the metric coincides with the one predicted by
Gross-Siebert \cite{GS} and Haase-Zharkov \cite{HZ}. 

\vspace{0.5cm}

\indent The LYZ construction of the metric comes down to looking for solutions
of the definite affine sphere equation \cite{SW}
\be
\label{LYZeq}
\psi_{z \bar{z}} + \frac{1}{2} e^{\psi} + |U|^{2} e^{-2\psi} = 0, \qquad
U_{\bar{z}}=0,
\ee
where $\psi$ and $U$ are real and complex functions respectively on
an open set in $\C.$  LYZ set $U= z^{-2}$ to account for the
singularity of the metric they considered.  They then proved the
existence of the radially symmetric solution $\psi$ of
(\ref{LYZeq}) with a prescribed behaviour near the
singularity $z=0,$ and established the existence of the global solution to the
coordinate-independent version of (\ref{LYZeq}) on $S^2$ minus three
points.

In this paper, we study the integrability of equation (\ref{LYZeq}).
We show that the affine sphere equation and a closely
related equation called the Tzitz\'eica equation arise as reductions of 
anti--self--dual Yang--Mills (ASDYM) system by two translations,
 and hence it admits a twistor interpretation.
Moreover, the ODE characterising its radial solutions gives rise to an
isomonodromy problem described by the Painlev\'e III ODE. The two-dimensional group of translations
reduces the Euclidean ASDYM equations  to the Hitchin equations
\cite{Hitchin_system} and Theorem \ref{Main_thm_LYZ} below
gives an invariant characterisation of (\ref{LYZeq})
as a special case of the $SU(2, 1)$ Hitchin equations.

Let $A$ be an $\mathfrak{su}(2,1)$ valued connection on a rank 3 complex vector
bundle $E\rightarrow \C$ with the
curvature $F_A=dA+A\wedge A$
and let $\Phi$ be a one-form with values in adj$(E)$. Choose a local trivialisation
of $E$ and set \[
A=A_z dz+ (A_z)^* d\bar{z},\quad
\Phi=Q d\bar{z}, \quad D = d +A,
\]
where $m^*:=-\eta^{-1} \bar{m}^t\eta$
with $\eta=\mbox{diag}(1, 1, -1)$, so that $\Phi^*=Q^* dz$.
\begin{theo} \label{Main_thm_LYZ}
The Hitchin equations
\be
\label{hitchin_eq_theorem}
F_A-\Phi\wedge \Phi^*-\Phi^*\wedge\Phi=0, \qquad
D\Phi=0
\ee
hold with
\be
\label{LYZ_Hichin}
A_z= \left( \begin{array}{ccc}
    0 & \frac{1}{\sqrt{2}}e^{\frac{\psi}{2}} & 0 \\
    0 & -\frac{1}{2}\psi_z & -Ue^{-\psi} \\
    0 & 0 & \frac{1}{2}\psi_z
  \end{array}
  \right), 
\qquad Q=
\left(
\begin{array}{ccc}
    0 & 0 & \frac{1}{\sqrt{2}}e^{\frac{\psi}{2}} \\
    0 & 0 & 0 \\
    0  & 0 & 0
  \end{array}
  \right)
\ee
if the functions $(\psi, U)$ satisfy the affine sphere equation
{\em(\ref{LYZeq})}. 

Conversely, any solution
to the $SU(2, 1)$ Hitchin equations such that
\begin{enumerate}
\item $Q$ has minimal polynomial $t^2$ 
and ${\mbox Tr}(Q Q^*) \ne 0,$ 
\item 
${\mbox Tr} \left( (D_z Q^*)^2 \right) = 0,\quad 
{\mbox Tr}\left( (D_z Q^*)^2(D_{\bar z}Q)^2\right) \ne 0,$
\item ${\mbox Tr} [(QQ^*)^4 -
 (Q^*Q)^2(D_zQ^*)(D_{\bar z}Q) +
Q^*Q(D_zQ^*)QQ^*(D_{\bar z}Q)]=0$
\end{enumerate} 
is equivalent to {\em(\ref{LYZ_Hichin})} by gauge and coordinate transformations.
\end{theo}

The connection between solutions to the affine sphere equation (\ref{LYZeq}) and the Calabi--Yau metric
(\ref{g_omega}) in six dimensions has not been made explicit in \cite{LYZ}. The Lax representation of 
(\ref{LYZeq}) will be used to prove the following
\begin{prop} \label{CYprop}
Given a semi-flat Calabi--Yau metric {\em {(\ref{g_omega})}}, where $\phi(x)$
satisfies the Monge--Amp\'ere equation {\em {(\ref{hessian})}}, and 
$\phi(cx)=c^2 \phi(x)$ where $c$ is a non--zero constant,
there exist complex coordinates $\{ z, w, \xi \}$ such that the metric
$g$ and the K\"ahler form $\om$ can be written as
\begin{eqnarray}
g &=& e_1 \b e_1 +  e_2 \b e_2 +  e_3 \b e_3,  \label{semiflatcymetric} \\
\om &=& \frac{i}{2} \left( e_1 \wedge \b e_1 + e_2 \wedge \b e_2 + e_3
\wedge \b e_3 \right), \nonumber
\end{eqnarray}
where
\begin{eqnarray}
e_1 &=& dw -\frac{i}{2}e^\psi(\b \xi dz + \xi d\b z), \nonumber \\
\label{basis}  e_2 &=& \frac{e^{\psi/2}}{\sqrt{2}} \left( \; (w+i\xi\psi_z)dz + i(d\xi +
e^{-\psi}\b U \b \xi d\b z) \;\right) \\
e_3 &=& \frac{e^{\psi/2}}{\sqrt{2}} \left(\;  i(d\b \xi +
e^{-\psi}U\xi dz) + (w+i\b \xi\psi_{\b z}) d\b z \; \right), \nonumber
\end{eqnarray}
and $\psi(z, \b z),$ $U(z)$ are real and complex functions
respectively defined on an open set in $\C$ which satisfy the 
affine sphere equation {\em {(\ref{LYZeq})}}.
\end{prop}  
The Hitchin equations (\ref{hitchin_eq_theorem}) are integrable as they arise from ASDYM
and their solutions can be described by holomorphic twistor data. Therefore any ODE arising as 
reduction of (\ref{LYZeq}) by another symmetry 
must be of Painlev\'e type in agreement with an integrable
dogma \cite{ARS, MW, D2009}.

If $U=z^n, \; n \in \Z,$ the equation (\ref{LYZeq})
admits rotational symmetry
\be \label{rotation} z\rightarrow e^{ic}z, \quad c \in \R. \ee
Therefore one can consider the group invariant solutions $\psi$
and look for the ODE characterising such reduction.
For concreteness, let us consider $U= z^{-2}$ following LYZ. 

\begin{prop}
\label{painleve_prop}
Solutions to {\em(\ref{LYZeq})} with
$U=z^{-2}$ invariant under
a group of rotations {\em(\ref{rotation})} are
of the form
\[
\psi(z, \bar{z})=\log{H(s)}-3\log{(s)},\qquad s=|z|^{1/2},
\]
where $H$ satisfies 
\[
H_{ss}=\frac{{(H_s)}^2}{H}-\frac{H_s}{s}-\frac{8H^2}{s} -\frac{16}{H}
\]
which is the Painlev\'e III equation with parameters $(-8, 0, 0, -16)$.
\end{prop}

In the next section we follow Leung \cite{Leung} and review the semi--flat
Calabi-Yau manifolds.  Then, in section 3 we summarise the results about affine spheres
which are used in the LYZ construction \cite{LYZ}.  In
section 4 we prove Theorem \ref{Main_thm_LYZ} and
give a gauge invariant characterisation of the definite
affine sphere equation
and the closely related Tzitz\'eica equation as symmetry reductions of the 
anti--self--dual Yang--Mills equations. As a byproduct, in section 5 we shall obtain a
characterisation of a reduction of the Hitchin equations to the $\Z_3$ two dimensional Toda chain. In section 6 we discuss other possible gauge inequivalent reductions of
the ASDYM equations to the affine sphere equation and the Tzitz\'eica equation.
In section 7  we  give a proof of Proposition \ref{CYprop} and
recover the toric
Calabi--Yau metric in terms of the solutions of the affine sphere equation.
Finaly in  section 8 we establish Proposition \ref{painleve_prop} and
demonstrate
that the existence theorem for Hessian metrics with prescribed
monodromy
comes down to the study of the Painlev\'e III equation with special
values of parameters, and obtain the corresponding $3 \times 3$
isomonodromic Lax pair.  
\section{Semi--Flat Calabi--Yau manifolds and the SYZ conjecture}
\setcounter{equation}{0}
Let $z^j=x^j+i y^j$ be holomorphic coordinates on a Calabi--Yau three-fold 
$X$, and let $\phi(z^j, \bar{z}^j)$ be the K\"ahler potential such that
$\omega =i \partial \overline{\partial}\phi$. The Ricci--flat condition for 
the corresponding Riemannian metric is 
\[
\Omega\wedge\overline{\Omega}=\omega^3,
\]
where $\Omega=d z^1\wedge d z^2\wedge d z^3$ is the holomorphic 
three-form on $X$. 

Now let us consider the $T^3$ invariant case.  Assume that the potential $\phi$ is invariant under translations in
the imaginary directions $y^j$.  In this case the Riemannian metric
and the K\"ahler form  are given by
(\ref{g_omega}) where 
\[
\phi_{jk}:=\frac{\p^2\phi}{\p x^j \p x ^k}
\]
and the Ricci--flat condition reduces to the real Monge--Amp\'ere equation 
(\ref{hessian}) for $\phi=\phi(x^1, x^2, x^3)$.

We shall regard the $x^j$ as local coordinates in an open set 
$B\subset \R^3$.
The freedom in choosing the coordinates $x^j$ without changing the equation
(\ref{hessian}) is given by affine transformations
${\bf x}\rightarrow {M\bf x}+{\bf b}$, where $M\in SL(3, \R)$, and 
${\bf b}$ is a vector. The affine transformations induce the change
in the potential $\phi\longrightarrow (\mbox{det}{M})^2\phi$, thus
$\phi$ should be regarded as a section of the second power of the 
real determinant line bundle over $B$. Conversely, given a three real
dimensional affine manifold $B$ with 
a metric of Hessian type (\ref{Hessian_met})
where $\phi$ satisfies the Hessian condition (\ref{hessian}) one can construct
the Calabi--Yau metric on $X=TB$ by (\ref{g_omega}). We then compactify
the fibres quotienting them by a lattice thus producing a $T^3$ 
invariant Calabi--Yau 
structure on the total space of a toric fibration
$\pi:X\longrightarrow B$. 

We are now ready to formulate the SYZ conjecture.
If $X, Y$ are mirror Calabi--Yau manifolds (see \cite{mirror} for a discussion
of what it means) then there exists a compact real 
three-manifold $B$ such that
\begin{itemize}
\item
$
\pi:X\longrightarrow B, \quad \rho:Y\longrightarrow B
$
are special Lagrangian fibrations by tori (the fibres can be singular
at some points of $B$).
\item The fibres of $\pi$ and $\rho$ are dual tori. 
\end{itemize} 
The second condition only makes sense for flat tori, therefore
the conjecture holds in the {\em large complex structure limit}, where
the volume of the fibres is small in comparison to the volume of the 
base space and the metric on the fibres is approximately flat.

To understand  the large complex structure limit consider a one parameter
family of complex structures $J(t)$ given by the holomorphic coordinates
\[
z^j(t)=t^{-1} x^j+i y^j,
\]
and the corresponding  Calabi--Yau metrics rescaled by $t^2$
\[
g(t)=\phi_{ij}(d x^j d x^k+t^2 d y^j d y^k).
\]
Thus we get a one parameter family of special Lagrangian fibrations.
In a limit $t\longrightarrow 0$ the Gromov--Hausdorff limit of 
metric $g(t)$ is the Hessian metric (\ref{Hessian_met}) on $B$, and
the size of the fibres shrinks to zero. The SYZ conjecture
predicts that such a limit exists for any Calabi--Yau metric
on a (not necessarily $T^3$ symmetric) toric special Lagrangian
fibration.

\section{Affine geometry and Hessian metrics}
\setcounter{equation}{0}
The Hessian equation (\ref{hessian}) is known not
to be integrable, at least in the sense of the
hydrodynamic reductions \cite{Fer}. Its homogeneous
solutions are however characterised by an integrable PDE. 
We shall carry over the homogeneity analysis
for a general Hessian metric in $(n+1)$ dimensions,
and then restrict our attention to $n=2$ where
there is a direct connection with the semi--flat
CY manifolds on one side and integrability on the other.

The following proposition follows from combining results of Calabi \cite{Calabi} 
and Baues-Cort\'es \cite{B-C}
about parabolic and elliptic affine spheres.
Here, we give a direct elementary proof not based on affine
differential geometry. It has certain advantages
as it exhibits explicit coordinate transformations between solutions
to various forms of homogeneous Hessian equations.
\begin{prop}
\label{prop_hessian}
Let $\phi=\phi(x^i)$ be a solution to the Hessian equation
{\em(\ref{hessian})} on an open ball $B\subset
\R^{n+1}$
such that $\phi(c x)=c^2 \phi(x)$ for any non-zero constant $c.$
Then there exists a local coordinate
system $(p_1, \dots, p_n, r)$ on $B$
such that the metric {\em(\ref{Hessian_met})} is
\be
\label{h_legendre}
g_B=dr^2+ r^2\frac{1}{w}\Big(\frac{\p^2 w}{\p p_{\alpha}\p p_{\beta}}
\Big)dp_{\alpha} dp_{\beta}, \qquad \alpha, \beta=1, \dots, n,\,
\ee
where $w=w(p_\alpha)$ satisfies 
\be
\label{eq_legendre}
\det{\Big(\frac{\p^2 w}{\p p_{\alpha}\p p_{\beta}}\Big)}=\frac{1}{w^{n+2}}.
\ee
\end{prop}
{\bf Proof.}
Consider the Hessian metric  (\ref{Hessian_met})
with $\phi$ homogeneous of degree 2. 
Therefore $V=x^i \p/\p x^i$ is a homothety with ${\cal L}_V g_B=2g_B.$
Locally there exists a function $r:B \lra \R$ such that $V= r \p/\p r$ and
\[
g_B= \gamma(dr+r\alpha)^2+r^2 h
\]
where $h, \alpha, \gamma $ are a metric, a one--form and a function
respectively  on the space of orbits of $V.$ The relation
$\p_i(x^j\phi_j)=2\phi_i$ gives
\[
g_B(V, ...)=x^i\phi_{ij}dx^j=d\phi.
\]
Thus $d(\gamma(dr+r\alpha))=0$ and we can redefine $r$ to set 
$\alpha=0$ and $\gamma=1$. We also note that
$|V|^2=x^ix^j\phi_{ij}=2\phi,$ and recognise 
$g_B$ as a cone over $h$ 
\be \label{r}
g_B=dr^2+r^2h, \quad \phi=\frac{r^2}{2}.
\ee
Now let us consider the surface $r=1$ given by a graph in $\R^{n+1}$
\[ (\t x^1, \dots, \t x^n) \longmapsto (\t x^1, \dots,\t x^{n}, v(\t x^\a)), \]
where $\t x^\a, \alpha=1, \dots, n,$ parametrise the surface. We shall
show that its induced metric $h$ is given by
\be
\label{Monge_metric_1}
h=\frac{\p_\alpha\p_\beta v}{\t x^\gamma\p_\gamma v-v}
d\t x^{\alpha} d\t x^{\beta},
\ee
where $\p_\a := \p/\p \t x^\a.$
To prove it, restrict the function $\phi$ to the surface $r=1.$  This gives an identity
$\phi(\t x^\alpha, v(\t x^\alpha))=1/2$. 
We differentiate this identity implicitly  with respect to $\t x^\alpha$ and express
the first and second derivatives of $\phi$ in terms of the  derivatives of $v$
\begin{eqnarray*}
0&=&\p_\alpha\phi+\p_{n+1}\phi\,\p_\alpha v,\\ 
0&=&\p_\alpha\p_\beta\phi+\p_\alpha\p_{n+1}\phi\,\p_\beta v
+\p_\beta\p_{n+1}\phi\,\p_\alpha v+\p_{n+1}^2 \phi\,
\p_\alpha v\p_\beta v+\p_{n+1} \phi\,
\p_\alpha\p_\beta v, \\
2\phi&=&\t x^\alpha\p_\alpha\phi+v\p_{n+1}\phi=1,
\end{eqnarray*}
where the last relation is just the homogeneity condition
restricted to the hypersurface $\phi=1/2$. Substituting all
that to $g_B$ gives (\ref{Monge_metric_1}).

Now if the function $\phi$ in the Hessian metric
$g_B$ satisfies the Hessian condition (\ref{hessian})
then $v$ satisfies
\be
\label{u_equation}
\det{
\frac{\p^2 v}{\p \t x^\alpha \p \t x^{\beta}}}=
(\t x^\alpha \p_\alpha v-v)^{n+2}.
\ee
To see it, let us write the coordinates $x^i$ on $\R^{n+1}$
 as $(x^1, \dots,
x^n, x^{n+1}) = (r\t x^1, \dots,r\t x^{n}, rv(\t x^\a)),$ that is,
regard $\R^{n+1}$ as the cone over the $r=1$ surface.  
Now consider the invariant volume element
\be
\label{vol_identity}
\sqrt{|g_B|}\,dx^1\wedge\ldots\wedge dx^n \wedge dx^{n+1}=
\sqrt{|\tilde{g}_B|}\,d\tilde{x}^1\wedge\ldots \wedge d\t x^n \wedge dr
\ee
where $|g_B|$ is the absolute value of the determinant of Hessian
metric (\ref{Hessian_met}) written in the  coordinates $x^i$
and $\tilde{g}_B$ is the same metric expressed in the basis $\{d\t
x^\a, dr\}.$
We contract both sides of $(\ref{vol_identity})$ with 
$V$. On the LHS of  (\ref{vol_identity}) we use the form 
$V= x^i \p/\p x^i$ and on the RHS  use  $V= r \p/\p r.$ 
We now set $r=1$ 
and impose the Hessian equation (\ref{hessian}), $\det{g_B}= \det \phi_{ij}=1.$ 
This yields
\[
v-\t x^{\alpha}\p_{\alpha} v=\sqrt{|\tilde{g}_B|}.
\]
On the surface $r=1,$ one has $\det{\tilde{g}_B}=\det{h}$ 
where $h$ is given by (\ref{Monge_metric_1}). Substituting this in the
above formula and taking squares of both sides yields 
(\ref{u_equation}).  Note\footnote{If we started with
$\det \phi_{ij}=-1,$ which implies $\det h < 0,$ the analogous
  argument would lead to  $\det{
\frac{\p^2 v}{\p \t x^\alpha \p \t x^{\beta}}}=
-(\t x^\alpha \p_\alpha v-v)^{n+2}.$} that we have taken $\det h >0$ from the
assumption that $\det g_B = \det \phi_{jk} =1.$

To obtain the statement in the proposition, perform  a Legendre transform
\[
p_\alpha=\frac{\p v}{\p \t x^{\alpha}}, \quad
w(p_\alpha)=\t x^{\alpha}\frac{\p v}{\p \t x^{\alpha}}-v,
\quad \t x^{\alpha}=\frac{\p w}{\p p_{\alpha}}.
\]
Using $dp_\alpha=\p_\alpha\p_\beta v\,d\t x^\beta$
yields 
\be
\label{Monge_metric_3}
h=\frac{1}{w}\frac{\p^2 w}{\p p_{\alpha}\p p_{\beta}} dp_{\alpha} dp_{\beta}
\ee
and
\[
{\frac{\p^2 w}{\p p_{\alpha}\p p_{\beta}}}=
\Big({\frac{\p^2 v}{\p \t x^{\alpha}\p \t x^{\beta}}}\Big)^{-1},
\]
which implies (\ref{h_legendre}) and (\ref{eq_legendre}).
\koniec
Now, let us consider a hypersurface $\Sigma$ immersed in
$\R^{n+1}$ with the flat metric $\delta_{jk}\; d x^j dx^k,$ given by a graph
\be \label{graph}
{\bf r}=(\t x^1, \dots, \t x^n, v(\t x^1, \dots, \t x^n)).
\ee
The first
and second fundamental forms on $\Sigma$ are given by
\begin{eqnarray*}
h_I&=&d{\bf r}\cdot d{\bf r}=(\delta_{\alpha\beta}
+\p_\alpha v\p_\beta v)d\t x^\alpha d\t x^\beta,\\
h_{II}&=&-d{\bf r} \cdot d{\bf n}=\frac{1}{\sqrt{1+(\p_1 v)^2+\dots
    +(\p_n v)^2}}\, \frac{\p^2 v}{\p \t x^\alpha\p \t x^\beta} d\t
x^\alpha d\t x^\beta,
\end{eqnarray*}
where ${\bf n}$ is the unit normal to $\Sigma$.
Tzitz\'eica \cite{Tzitzeica1,Tzitzeica2} has studied surfaces $\Sigma$ in $\R^3$
for which the ratio of the Gaussian curvature ${\cal K}$ to the 
fourth power of a distance from a tangent plane to some fixed point is
a constant.  If ${\cal K} \ne 0,$ we can always rescale the coordinates to set this constant
to $+1$ or $-1$ depending on the sign of the Gaussian curvature.
We shall call this the Tzitz\'eica condition.  The generalisation of the
Tzitz\'eica condition to
hypersurfaces in $\R^{n+1}$ is given by
\[ {\cal K} = \pm {\cal D}^{n+2}, \]
where ${\cal D}= {\bf r}\cdot{\bf n}$ is the same as the distance up to sign. 
In the adapted coordinates, ${\cal D}$ 
and the Gaussian curvature ${\cal K}$ are given by
\begin{eqnarray*}
{\cal D}&=&\frac{v-\t x^{\alpha}\p_\alpha v}{\sqrt{1+(\p_1 v)^2+\dots +(\p_n v)^2}},\\
{\cal K}&=&\frac{1}{(\sqrt{1+(\p_1 v)^2+\dots +(\p_n v)^2})^{n+2}}\,
\det{\Big(\frac{\p^2 v}{\p \t x^\alpha\p \t x^\beta}\Big)}.
\end{eqnarray*}
It follows that the Tzitz\'eica condition holds if and only if $v$ satisfies
\be \label{Tcond} \det{
\frac{\p^2 v}{\p \t x^\alpha \p \t x^{\beta}}}= \pm
(v-\t x^\alpha \p_\alpha v)^{n+2},\ee
where plus and minus signs correspond to positive and negative Gaussian curvature
respectively.  

It is well known in affine differential geometry that an immersed
hypersurface $\Sigma$ in $\R^{n+1}$ is an affine hypersphere with the origin as
its centre if and only if the Tzitz\'eica condition (\ref{Tcond}) holds
\cite{NS}.  It turns out that the metric
(\ref{Monge_metric_1}), with $v$ satisfying (\ref{u_equation}),
is the same as the Blaschke metric (or affine
metric) of a proper affine hypersphere.  The Blaschke metric is
conformally related to the second fundamental form, and is defined as follows.  
Let ${\bf N}$ denote the transversal vector field of the surface
$\Sigma$ such that the unit normal ${\bf n}$ is given by
${\bf n} = \frac{\bf N}{|{\bf N}|},$
i.e. ${\bf N}=\nabla(\t x^{n+1}-v(\t x^1,...,\t x^n)).$ Consider a
bilinear form 
\[\h h = -d{\bf r}\cdot d{\bf N} = |{\bf N}| \; h_{II}. \]
The Blaschke metric is then given by
\be \label{Blaschke} h := |\det \h h|^{-\frac{1}{n+2}} \; \h h.\ee
Therefore, for the surface $\Sigma$ given by the graph (\ref{graph}),
we have
\[
h= \Big| \det{\frac{\p^2 v}{\p \t x^\alpha\p \t x^\beta}}\Big|^{-\frac{1}{n+2}}
\frac{\p^2 v}{\p \t x^\alpha\p \t x^\beta} d\t x^\alpha d\t x^{\beta},
\]
which coincides with the metric (\ref{Monge_metric_1}) if equation (\ref{u_equation}) holds. 

In affine differential geometry, it is also known \cite{Calabi} that a Hessian metric 
(\ref{Hessian_met}) which satisfies $\det{\phi_{ij}}=1$ is a
parabolic (improper) affine hypersphere metric. 
We have demonstrated that Hessian equation (\ref{hessian}) on $\phi$
implies (\ref{u_equation}) on $v$.  Therefore, this is in agreement with a result of
Baues and Cort\'ez \cite{B-C} that a parabolic affine hypersphere metric
which admits a homothety ${\cal L}_V g_B=2g_B$ is the metric cone over a
proper affine hypersphere.  

Let us now restrict our attention to $n=2,$ and consider the metric
$h$ (\ref{Monge_metric_1}).
For $n=2,$ $\det h > 0$ implies that $h$ is a definite metric.  
In the context of the Calabi--Yau manifolds, the metric $g_B$ is
Riemannian, hence one is interested
in positive--definite $h.$  Baues and Cort\'es \cite{B-C} have shown
that in such case $h$ is the Blaschke metric of a definite elliptic
affine sphere, with affine mean curvature 1.  Since $h$ is positive definite
we can adopt isothermal coordinates
for the affine  metric (which are asymptotic coordinates for the second
fundamental form $h_{II}$) and write it as
\be \label{affineh}
h=e^{\psi}dz d\overline{z},
\ee
for some real valued function $\psi=\psi(z, \b z).$
In this form, Simon and Wang \cite{SW} proved that the structure
equations\footnote{The usual affine immersion in
  $\R^{n+1}$ only assumes a flat connection $D$ and a parallel
  volume element on $\R^{n+1},$ but not an ambient metric.  In
  particular, the structure equations
  of a Blaschke hypersurface immersion $f: (\Sigma, \nabla) \lra
  (\R^{n+1}, D)$ are given by
\be \label{Gauss} D_X f_* (Y) = f_*(\nabla_X Y) + h(X,Y)\xi,  \ee
\be \label{Weingarten} D_X \xi = -f_*(SX), \ee
where $\nabla$ is an affine connection on $\Sigma,$ $X, Y \in
T\Sigma, \; \xi$ is a transversal vector field chosen uniquely up to
sign to satisfy certain properties, called the affine normal field,
and $h$ is the Blaschke metric defined by (\ref{Gauss}).
This definition turns out to be
  equivalent to (\ref{Blaschke}) if one were to use the Euclidean
  metric on $\R^{n+1}.$  The operator $S: T\Sigma \lra T\Sigma$ is
  called the affine shape operator and $H= \frac{1}{n}\mbox{Tr}(S)$
  the affine mean curvature.  A proper affine sphere is defined to be a
  Blaschke hypersurface with $S=HI,$ $I$ being the identity metric.
 Another affine invariant quantity
  is a totally symmetric tensor called the cubic form $\h C$ and is
  defined by
\[\h C (X,Y,Z) = h(C(X,Y),Z),\]
where $C$ is the difference tensor $C= \h \nabla - \nabla$ and $\h
\nabla$ is the Levi-Civita connection of $h.$  Consider $h$ as in
(\ref{affineh}) and let $C^i_{jk}, \; i,j,k \in \{1, \b 1 \}$ be the components of $C$
in the basis $e^1=dz, e^{\b 1}=d\b z.$  Then it can be shown that the only
nonvanishing components of $C$ are $C^{\b 1}_{11}$ and $C^{1}_{\b 1 \b
  1}= \overline{C^{\b 1}_{11}},$ and the function $U$ in (\ref{LYZeq})
is defined by $U = C^{\b 1}_{11} e^\psi.$  It follows that the cubic
form is $\hat C =  U dz^3 + \b U d \b z^3.$  See \cite{Calabi, NS, SW, loftin_review}
for details.} 
of definite affine sphere imply that $\psi$ necessarily
satisfies the equation (\ref{LYZeq})
\[
\psi_{z \bar{z}} + \frac{1}{2} e^{\psi} + |U|^{2} e^{-2\psi} = 0, \qquad
U_{\bar{z}}=0,
\]
where $U dz^3$ is the holomorphic cubic differential.

Conversely, given a solution of (\ref{LYZeq})
one can construct an affine sphere with $h = e^{\psi}dz d\overline{z}$
as its Blaschke metric.  We should note here that if the holomorphic cubic
differential $U(z)dz^3$ is non-zero, we can choose the
isothermal coordinates such that $U=1.$  For example, defining $\xi=\xi(z)$ by
$ d\xi=2^{-1/3}U^{1/3}dz$ transforms (\ref{LYZeq}) into
\be \label{ellipticLYZ}
\hat{\psi}_{\xi\overline{\xi}}+e^{\hat{\psi}}+e^{-2\hat{\psi}}=0,
\ee
where
\[
\hat{\psi}={\psi}-\frac{1}{3}\log{U}-\frac{1}{3}
\log{\bar{U}}-\frac{1}{3}\log{2}.
\]
We will make use of such coordinate transformation in section 4\footnote{ We note that the analytic continuation
\[
\hat{\psi}_{\xi\overline{\xi}}+e^{\hat{\psi}}-e^{-2\hat{\psi}}=0
\]
of equation (\ref{ellipticLYZ}) was used by  McIntosh \cite{mcintosh} to describe minimal Lagrangian immersions in
$\CP^2$ and special Lagrangian cones in $\C^3$.}. 
\vspace{0.5cm}
\indent Loftin, Yau and Zaslow \cite{LYZ} proved the existence of a semi--flat Calabi--Yau
metric (\ref{g_omega}) with the base metric $g_B$ as the metric cone
over an elliptic affine sphere
\be \label{gB}
g_B=\phi_{ij}dx^idx^j=dr^2+r^2 e^\psi dz d\overline{z},
\ee
with the prescribed singularity, by proving the existence of a radially symmetric solution $\psi$ of
(\ref{LYZeq}) for $U(z)=z^{-2}$ and the corresponding global solution on $S^2$ minus three
points. 

Motivated by this work, we are interested in the integrability of the
definite affine sphere
equation (\ref{LYZeq}). The affine sphere equation is closely related to a well
known integrable equation, namely the Tzitz\'eica equation
\be \label{Teqn} 
u_{xy}=e^u- e^{-2u}.
\ee

In the context of affine spheres, the Tzitz\'eica equation arises if $\det h < 0.$  
By writing the metric in isothermal
coordinates as $h=2e^u\;dxdy$ and considering the structure equations,
Simon and Wang \cite{SW} also show that $h$ is the Blaschke metric of the indefinite
affine sphere (with negative affine mean curvature) if and only if $u$
satisfies $ u_{xy}=e^u - r(x)b(y) e^{-2u},$ where $r(x), b(y)$ are
arbitrary non-vanishing functions of one variable, which can be
normalised by rescaling the isothermal
coordinates.  Thus, we obtain
\be \label{definitesph} u_{xy}=e^u-\epsilon e^{-2u},\ee
where $\epsilon = \pm 1.$  The equation with $\epsilon=1,$ (\ref{Teqn}), was first
derived in \cite{Tzitzeica1,Tzitzeica2}
for the Tzitz\'eica surface in $\R^3$ with negative Gaussian curvature
${{\cal K}=-{\cal D}^4,}$ where the indefinite second fundamental form
is written in asymptotic
  coordinates as $h_{II} = 2e^u{\cal D}\;dxdy.$  

The difference between the two equations (\ref{Teqn}) 
and (\ref{LYZeq}) lies
in the relative sign of the two exponential terms on the RHS.
For the Tzitz\'eica equation $u=0$ is a solution and
other solutions may be constructed using Darboux and B\"acklund
transformations, for example see \cite{BSS}. The
definite affine sphere equation does not seem to have such obvious solutions.
However, Calabi \cite{Calabi} has shown that an elliptic
affine hypersphere with complete Blaschke metric is an ellipsoid.
This is in agreement with the fact that (\ref{LYZeq}) admits solutions
in term of elliptic functions, which can be found by making an ansatz
$\psi(z, \b z)=f(z +\b z)$ in (\ref{ellipticLYZ}).

\section{Reduction of ASDYM}
\setcounter{equation}{0}
It was shown in \cite{Dunajski} that the 
Tzitz\'eica equation (\ref{Teqn}) can be obtained from a special
ansatz to the anti--self--dual Yang--Mills in $\R^{2,2}$ with gauge
group $SL(3, \R).$ In this section, we
shall give a gauge and coordinate invariant characterisation of
 the Tzitz\'eica equation and the definite affine sphere equation as 
different real forms of a
reduction of ASDYM on $\C^4$ with gauge
group $SL(3, \C),$ via
the holomorphic Hitchin equations on $\C^2.$

\subsection{Holomorphic Tzitz\'eica equation}
Consider a holomorphic metric and volume element on  $\C^{4}$
$$
ds^{2} = 2(dz \ d\tilde{z} - dw \ d\tilde{w}),\quad
\nu = dw \wedge d\tilde{w} \wedge dz \wedge d\tilde{z}.
$$
Let ${\cal A}=A_zdz+A_wdw+A_{\tilde{z}}d\tilde{z}+A_{\tilde{w}}d\tilde{w}$ 
be a Lie algebra valued connection on a vector bundle $E\rightarrow \C^4$.
The anti--self--dual Yang--Mills equations are
given by
\[
F_{zw}=0, \quad 
F_{z \tilde{z}} - F_{w \tilde{w}} =0, \quad F_{\tilde{z} \tilde{w}} =0.
\]
These equations arise from a Lax pair
\be
\label{ASDYM_Lax}
[D_z+\lambda D_{\tilde{w}}, D_w+\lambda D_{\tilde{z}}]=0,
\ee
where  $D_z=\p_z+A_z,$ etc, are covariant derivatives,
$F_{z\tilde{z}}=[D_z, D_{\tilde{z}}],$
and (\ref{ASDYM_Lax}) is required to hold for any
value of the spectral parameter $\lambda$.

Choose a gauge group to be $SL(3, \C)$ 
and assume that ${\cal A}$ is 
invariant under the action of two dimensional group of  translations
$\C^2$ such that the metric restricted to the planes spanned by
the generators of the group is non-degenerate. 
Let $X_1, X_2$ be the generators of the group, then the Higgs fields 
\[
P=X_1\hook A,\quad Q= X_2\hook A
\]
belong to the adjoint representation.
We can always choose the coordinates so that the group is
generated by the two null vectors
$X_1=\p/\p \tilde{w}$ and $X_2=\p/\p w$.
The ASDYM system reduces
to the holomorphic form   of the Hitchin equations \cite{Hitchin_system}
\begin{subequations}
\begin{align}
\label{holHitchin1}
D_z Q ={}& 0, \\
\label{holHitchin2}
D_{\t z} P ={}& 0, \\
\label{holHitchin3}
F_{z \t z}+[P, Q] ={}& 0,
\end{align}
\end{subequations}
where 
\[
F_{z\tilde{z}}=\p_z A_{\tilde{z}}-\p_{\tilde{z}} A_z+[A_z, A_{\tilde{z}}]
\] 
is a curvature of a holomorphic connection
$A=A_z dz +  A_{\tilde{z}} d{\tilde{z}}$  on $\C^2$.
The Hitchin equations are invariant under the gauge transformations
\be \label{gaugetransform}
A\rightarrow g^{-1} A g +g^{-1}\, dg, \qquad
P\rightarrow  g^{-1} P g, \quad
Q\rightarrow g^{-1} Q g
\ee 
and
later 
we shall also make use of the following coordinate freedom
\be 
\label{coordsym} z \lra \h z(z), \quad \t z \lra \h {\t z}(\t z). 
\ee

The Lax pair (\ref{ASDYM_Lax}) for the ASDYM 
 reduces to the following
Lax pair for the holomorphic Hitchin equations
\be
\label{lax_pair_H}
[D_z+\lambda P, Q+\lambda D_{\tilde{z}}]=0.
\ee
There are  several  gauge inequivalent ways to embed
the Tzitz\'eica equation (\ref{Teqn}) as a special case of the Hitchin equations.
The gauge used in \cite{Dunajski} is
\be
\label{holMD_Higgs}
A_{\tilde{w}}=P = \left(
\begin{array}{ccc}
0  & 0  & 1  \\
0 & 0  & 0  \\
0 & 0  & 0  
\end{array}
\right), \quad
A_w=Q = \left(
\begin{array}{ccc}
0  & 0  & 0  \\
0 & 0  & 0  \\
e^u & 0  & 0  
\end{array}
\right),
\ee
\be
\label{holMD_potential}
A_z = \left(
\begin{array}{ccc}
u_z  & 0  & 0  \\
1 & -u_z  & 0  \\
0 & 1  & 0  
\end{array}
\right), \quad
A_{\t z} = \left(
\begin{array}{ccc}
0  & e^{-2u}  & 0  \\
0 & 0  & e^u  \\
0 & 0  & 0  
\end{array}
\right),
\ee
where $u(z, \t z)$ is a complex valued function holomorphic
in $(z, \t z)$. With this ansatz the Hitchin equations yield
the holomorphic Tzitz\'eica equation
\be \label{holTeqn}
u_{z \t z}=e^u-e^{-2u}.
\ee
Choosing  the real form $SL(3, \R)$ of $SL(3, \C)$
and regarding $u=u(x, y)$ as a real function of
real coordinates $z=x, \tilde{z}=y$
reduces  (\ref{holTeqn}) to  (\ref{Teqn}). 

On the other hand, performing the coordinate transformation
\[ 
d {\h z} = \left(\frac{U(z)}{2} \right)^{-\frac{1}{3}} dz, \quad 
d {\h {\t z}} = \left(\frac{\t U(\t z)}{2} \right)^{-\frac{1}{3}} d{\t z}
\]
and setting  
\[ 
u = \psi(z, \tilde{z}) -\frac{1}{3}\log \left( \frac{U}{2} \right)  
-\frac{1}{3}\log \left( \frac{\t U}{2} \right)+ \log \left( -\frac{1}{2} \right) 
\]
for any branch of $\log \left( -\frac{1}{2} \right)$
puts (\ref{holTeqn}) in the form
\be
\label{holLYZeq}
\psi_{z \t z} + \frac{1}{2} e^{\psi} + U(z)\t U(\t z) e^{-2\psi} = 0,
\ee
where we have
dropped hats of the new variables.
Equation (\ref{holLYZeq}) then reduces to the affine sphere equation 
(\ref{LYZeq}) under the Euclidean reality
conditions $\tilde{z}=\bar{z}$ and reducing the gauge group to 
$SU(2, 1),$ which implies the constraint $\t U = \b U.$ 

\vspace{0.5cm}

\indent Now we shall establish a gauge invariant characterisation of
the ansatz (\ref{holMD_Higgs}), (\ref{holMD_potential}) in terms
of the gauge and Higgs fields of the Hitchin equations.
We will make use of the following lemma.  
\begin{lemma} \label{PQ}
Consider {\em 3} by {\em 3} complex matrices $P, Q$ such that 
\be
\label{conditions_1}
P^2=Q^2=0,\qquad Tr(PQ)=\omega \neq 0.
\ee
There exists a gauge transformation such that $P, Q$ are in the form
{\em (\ref{holMD_Higgs})} for some $u$.
\end{lemma}
{\bf Proof.} 
The conditions (\ref{conditions_1}) are invariant under the gauge 
transformations
\[
P\longrightarrow g^{-1} P g, \qquad Q\longrightarrow g^{-1} Q g.
\]
These conditions imply that the
nullities (dimensions of the kernels of the associated linear maps)
 satisfy $n(QP)<3$ and $n(P)=2$. Thus 
\[
\mbox{Ker}(QP)=\mbox{Ker}(P).
\]
Also rank$(QP)=1$ and Im$(QP)$ is contained in the one-dimensional image 
of $Q$,
therefore
\be
\label{imPQ}
\mbox{Im}(QP)=\mbox{Im}(Q).
\ee
Choose a Jordan  basis $({\bf v}, {\bf u}, {\bf w})$ of $\C^3$ such that
\be
\label{jordan_basis}
P({\bf w})={\bf v},\quad P({\bf v})=0,\quad P({\bf u})=0.
\ee
From (\ref{imPQ}) Im$(Q)=\mbox{span}\,(Q({\bf v}))$, thus 
$Q({\bf u})=a Q({\bf v}), Q({\bf w})=b Q({\bf v})$
for some $a, b$ so  that 
Ker$(Q)=\mbox{span}\, ({\bf u}-a {\bf v}, {\bf w}-b {\bf v})$.
Use the freedom in the basis (\ref{jordan_basis}) to set
\[
{\bf w'}={\bf w}-b {\bf v},\quad  {\bf u'}={\bf u}-a {\bf v},\quad 
{\bf v'} ={\bf v}.
\]
Now 
\begin{eqnarray*}
P({\bf w'})&=&{\bf v'}, \quad P({\bf v'})=0,\quad P({\bf u'})=0, \\
Q({\bf w'})&=&0,\quad Q({\bf u'})=0,\quad 
Q({\bf v'})= c {\bf u'}+ \om {\bf w'},
\end{eqnarray*}
where $\om \neq 0$ as $Tr(PQ)=\om \neq 0$. There is still freedom in 
(\ref{jordan_basis}):
\[
{\bf v''}={\bf v'},\quad {\bf u''}={\bf u'},\quad 
{\bf w''}={\bf w'}+ (c/\om){\bf u'} 
\]
so that, dropping primes,
\[
P({\bf w})={\bf v}, \quad P({\bf v})=0,\quad P({\bf u})=0,
\]
\[
Q({\bf w})=0,\quad Q({\bf u})=0, \quad Q({\bf v})= \om{\bf w}.
\]
Ordering the basis $({\bf v, u, w})$ yields the matrices in the desired form,
i.e. $P_{13}=1, Q_{31}=\om$, and all other components vanish.
The residual gauge freedom is 
\[
{\bf w}\rightarrow \alpha{\bf w},\quad {\bf v}\rightarrow \alpha{\bf v},\quad {\bf u}
\rightarrow \beta{\bf u} 
\]
and the change of
basis matrix gives the residual $GL(3, \C)$ 
gauge transformation. In the $SL(3, \C)$
case we set $\beta=\alpha^{-2}$.  The statement of the Lemma now follows by
setting $\om=e^u$.
\koniec 
We shall now give a set of necessary and sufficient conditions allowing
solutions of the Hitchin equations (\ref{holHitchin1}, b, c) to be transformed into 
(\ref{holMD_Higgs}), (\ref{holMD_potential}) by gauge and
coordinate symmetries.

\begin{prop} \label{gaugecond}
Let $(Q, P, A= A_zdz +A_{\t z}d{\t z})$ be a solution of the holomorphic Hitchin
equations {\em{(\ref{holHitchin1}, b, c})},
 with gauge group $SL(3,\C).$  
Then, $(Q, P, A_z, A_{\t z})$ can be transformed into
{\em{(\ref{holMD_Higgs}),(\ref{holMD_potential})}} by gauge symmetry and
coordinate symmetry {\em{(\ref{coordsym})}} if and only
if the following conditions hold:\\
${\bf (i)}$ $P$ and $Q$ have minimal polynomial $t^2,$ with ${\mbox
Tr}(PQ) \ne 0.$ \\
${\bf (ii)}$ ${\mbox Tr} \left( (D_zP)^2 \right) = 0 = {\mbox Tr} \left(
(D_{\t z}Q)^2 \right)$ and ${\mbox Tr}\left( (D_zP)^2(D_{\t z}Q)^2
\right) \ne 0.$
\\
${\bf (iii)}$ ${\mbox Tr}M=0,$ where \\
${} \qquad M = (PQ)^4 + (PQ)^2(D_zP)(D_{\t z}Q) -
PQ(D_zP)QP(D_{\t z}Q).$
\end{prop}
{\bf Proof.}
The proof of the necessary conditions is straightforward.
It can be shown by direct calculation that
(\ref{holMD_Higgs}),(\ref{holMD_potential})
 satisfy conditions ${\bf (i),(ii),(iii)}.$ The
three conditions are gauge
invariant by the cyclic property of the trace.  Under the coordinate transformation
(\ref{coordsym}), the 
connection $(A_z, A_{\t z})$ and the Higgs fields $(P,Q)$ transform as
\begin{eqnarray*}
 \h A_{\h z} &=& \left( \frac{d {\h z}}{d z} \right)^{-1}A_z,
   \quad \h A_{\h {\t z}} = \left( \frac{d {\h {\t
   z}}}{d {\t z}} \right)^{-1}A_{ \t z}, \\
\h Q &=& \left( \frac{d {\h {\t z}}}{d {\t z}} \right)^{-1}Q,
\quad \h P = \left( \frac{d {\h z}}{d z} \right)^{-1}P.
\end{eqnarray*}
Thus, using condition ${\bf (i)},$ the square of the covariant derivative
is given by
\[ (\h D_{\h z} \h P)^2 = \left( \frac{d {\h z}}{d z} \right)^{-4}
   (D_z P)^2 \]
and similarly for $(D_{\t z}Q)^2.$  Therefore, conditions ${\bf (i)}$ and
${\bf (ii)}$ are invariant under the coordinate transformation.  
A similar calculation shows that ${\bf (iii)}$ is also invariant under
(\ref{coordsym}).

Conversely, we shall now show that 
any solution to (\ref{holHitchin1}, b, c)  such that
all the conditions in Proposition
\ref{gaugecond} hold, can be gauge and coordinate transformed into
the form (\ref{holMD_Higgs}),(\ref{holMD_potential}).  

Firstly, by Lemma \ref{PQ}, condition ${\bf (i)}$ implies that
we can use gauge symmetry to put the Higgs fields $(Q,P)$ in the form (\ref{holMD_Higgs}).
The equations (\ref{holHitchin1}) and (\ref{holHitchin2}) imply that
$A_z, A_{\t z}$ are of the form 
\be \label{gaugecompat}
  A_z = \left( 
  \begin{array}{ccc}
    n & 0 & 0 \\
    r & u_z-2n & 0 \\
    m & t & n - u_z
  \end{array}
  \right), \quad
  A_{\tilde{z}} = \left( 
  \begin{array}{ccc}
    p & s & h \\
    0 & -2p & k \\
    0 & 0 & p
  \end{array}
  \right),
\ee
where $n,r,m,t,p,s,h,k$ are some functions of $(z, \t z).$  Note that  
we have also used the assumption that the fields are $\mathfrak{sl}(3,
\C)$ valued, hence traceless.  Next, to set the diagonal elements
of $(A_z, A_{\t z})$ to be as in (\ref{holMD_potential}), we consider 
the residual gauge freedom.   Lemma \ref{PQ} implies that the gauges
preserving $(Q,P)$ are given by
\be \label{gdiag}
 g(z,\t z) =  \left( 
  \begin{array}{ccc}
    a & 0 & 0 \\
    0 & \frac{1}{a^2} & 0 \\
    0 & 0 & a
  \end{array}
  \right)\ee
for an arbitrary function $a(z,\t z) \ne 0.$ Thus, using
(\ref{gaugetransform}), we have
\begin{eqnarray*}
  A_z &\longrightarrow& \left( 
  \begin{array}{ccc}
    n+\frac{a_z}{a} & 0 & 0 \\
    ra^3 & u_z-2n-2\frac{a_z}{a} & 0 \\
    m & \frac{t}{a^3} & n - u_z+ \frac{a_z}{a}
  \end{array}
  \right), \\
   A_{\t z} &\longrightarrow& \left( 
  \begin{array}{ccc}
    p+\frac{a_{\t z}}{a} & \frac{s}{a^3} & h \\
    0 & -2p-2\frac{a_{\t z}}{a} & ka^3 \\
    0 & 0 & p+\frac{a_{\t z}}{a}
  \end{array}
  \right).
\end{eqnarray*}
We choose $a(z,\t
z)$ such that
\[ (\ln a)_z = u_z - n, \quad {\mbox {and}} \quad (\ln a)_{\t z} =
-p.  \] 
This is allowed because the compatibility condition
\be \label{compatibility} \p_z p + \p_{z}\p_{\t z}u - \p_{\t z}n =0  \ee
holds automatically as a consequence of condition ${\bf(iii)}$.  To see it, note that
equation (\ref{holHitchin3}) implies 
\[ \p_z p + \p_{z}\p_{\t z}u - \p_{\t z}n +
mh + tk = e^u.     \]
Hence, condition (\ref{compatibility}) is equivalent to 
\[ mh + tk = e^u, \]
which holds by ${\bf (iii)}.$  

Note that at this point elements of $(A_z, A_{\t z})$ will be transformed, however,
for convenience we will label them with the same letters as in
(\ref{gaugecompat}).  Thus we have set $n=u_z$ and $p=0.$  We now proceed
to deal with $r, m, t, s, h, k.$  $\mbox{Tr}\left( (D_zP)^2(D_{\t z}Q)^2
\right) \ne 0$ in condition ${\bf (ii)}$ implies that
$r,t,s,k \ne 0,$ and 
\[
\mbox{Tr} \left( (D_zP)^2 \right) = 0 = \mbox{Tr} \left(
(D_{\t z}Q)^2 \right)\] 
gives
\[ m=0=h. \]
Hence (\ref{holHitchin3}) becomes
\begin{eqnarray*}
 u_{z \t z} + rs &=& e^u \\
 s_z + 2su_z &=& 0 \\
 r_{\t z} &=& 0 \\
 k_z -ku_z &=& 0 \\
 t_{\t z} &=& 0 \\
 tk &=& e^u.
\end{eqnarray*}
Since $r,t,s,k \ne 0,$ we can solve the above equations.  The last
three equations imply that $t$ is a constant, and thus can be set to $1$ by a
constant gauge transformation of the
form (\ref{gdiag}) with $a=t^{-1/3},$ and $s$ is determined to be of
the form $b(\t z)e^{-2u}.$   This results in
\[
P = \left(
\begin{array}{ccc}
0  & 0  & 1  \\
0 & 0  & 0  \\
0 & 0  & 0  
\end{array}
\right), \quad
Q = \left(
\begin{array}{ccc}
0  & 0  & 0  \\
0 & 0  & 0  \\
e^u & 0  & 0  
\end{array}
\right),
\]
\be \label{gMDansatz}
A_z = \left(
\begin{array}{ccc}
u_z  & 0  & 0  \\
r(z) & -u_z  & 0  \\
0 & 1  & 0  
\end{array}
\right), \quad
A_{\t z} = \left(
\begin{array}{ccc}
0  & b(\t z)e^{-2u}  & 0  \\
0 & 0  & e^u  \\
0 & 0  & 0  
\end{array}
\right).
\ee
Note that the gauge is now fixed.  To get to ansatz
(\ref{holMD_Higgs}),(\ref{holMD_potential}), we will now 
use the coordinate symmetry. Define $\h z, \h {\t z}$  such that
\[ d \h z = e^{j(z)} dz, \quad  d \h {\t z} = e^{l(\t z)} d \t z, \]
and set 
\[ \hat u := u - j(z) - l(\t z). \]
By choosing $j(z), l(\t z)$ such that $e^{3j(z)}=r(z)$ and $e^{3l(\t z)}=b(\t z),$
(\ref{gMDansatz}) becomes gauge equivalent to 
(\ref{holMD_Higgs}),(\ref{holMD_potential}) in the new variables 
$(\h z, \h {\t z}, \h u).$ 
The gauge transformation we need in the final step is given
by 
(\ref{gaugetransform}) with
\[ g(\h z, \h {\t z}) =  \left( 
  \begin{array}{ccc}
    e^{-j(z(\h z))} & 0 & 0 \\
    0 & e^{j(z(\h z))} & 0 \\
    0 & 0 & 1
  \end{array}
  \right).
 \]   \koniec
We note 
that substituting (\ref{gMDansatz}) to the Hitchin equations yields 
\be \label{gTeqn}
u_{z \t z} = e^u - r(z)b(\t z)e^{-2u}. 
\ee
Therefore, the change of coordinates can, roughly speaking, be regarded as
setting $r(z)$ and $b(\t z)$ to constants such that $r(z)b(\t z)=1.$

We shall now choose the Euclidean reality condition as
and select the real form $SU(2, 1)$ of $SL(3, \C)$ to
deduce Theorem \ref{Main_thm_LYZ} from the last Proposition.
{\bf Proof of Theorem \ref{Main_thm_LYZ}}.
 Consider the  ansatz (\ref{gMDansatz}) and equation
(\ref{gTeqn}).  By changing the dependent variable from $u$ to 
\[ \psi = u - \log \left( -\frac{1}{2} \right)  \]
for any branch of $\log \left( -\frac{1}{2} \right),$ equation (\ref{gTeqn}) becomes
\be \label{holLYZ} \psi_{z \t z} + \frac{1}{2}e^\psi + U(z)\t U(\t z)e^{-2 \psi} = 0, \ee
where $U(z)=2r(z), \; \t U(\t z) = 2b(\t z).$  Then, after an $SL(3,
\C)$ gauge transformation with
\[ g( z, {\t z}) =  \left( 
  \begin{array}{ccc}
    0 & 0 & -\sqrt{2}e^{-\frac{\psi}{2}} \\
    0 & \frac{1}{\sqrt{2}}e^{\frac{\psi}{2}} & 0 \\
    1 & 0 & 0
  \end{array}
  \right),
 \]
the ansatz (\ref{gMDansatz}) becomes
\begin{eqnarray} \label{MDLYZ}
 A_{w}&=&Q=\left(
\begin{array}{ccc}
    0 & 0 & \frac{1}{\sqrt{2}}e^{\frac{\psi}{2}} \\
    0 & 0 & 0 \\
    0 & 0 & 0
  \end{array}
  \right), \nonumber \\
  A_{\t w} &=& P= \left( 
  \begin{array}{ccc}
    0 & 0 & 0 \\
    0 & 0 & 0 \\
    -\frac{1}{\sqrt{2}}e^{\frac{\psi}{2}} & 0 & 0
  \end{array}
  \right),  \\  
A_{z}&=&\left( 
  \begin{array}{ccc}
    0 & \frac{1}{\sqrt{2}}e^{\frac{\psi}{2}} & 0 \\
    0 & -\frac{1}{2}\psi_z & -U(z)e^{-\psi} \\
    0 & 0 & \frac{1}{2}\psi_z
  \end{array}
  \right), \nonumber \\
  A_{\t z} &=& \left( 
  \begin{array}{ccc}
    0 & 0 & 0 \\
    -\frac{1}{\sqrt{2}}e^{\frac{\psi}{2}} & \frac{1}{2}\psi_{\t z}  & 0 \\
    0 & -\t U(\t z)e^{-\psi} & -\frac{1}{2}\psi_{\t z}
  \end{array}
  \right). \nonumber
\end{eqnarray}
Impose the Euclidean reality conditions
$\tilde{z} = \bar{z}$, $\tilde{w} = 
- \bar{w}$, resulting in a positive-definite metric on $\R^4$. The ASDYM 
equations with these reality conditions are 
\begin{eqnarray}
    F_{zw} &=& 0, \label{first} \\
    F_{z \bar{z}} + F_{w \bar{w}} &=& 0. \label{second}
\end{eqnarray}
Take the gauge group to be $SU(2,1)$. A matrix ${\cal M}$ is in the Lie 
algebra $\mathfrak{su}(2,1)$ if it is trace-free and satisfies
\be \label{liealgebra}
\bar{{\cal M}}^{t}=-\eta \; {\cal M}\; \eta^{-1},
\ee
where
$$
\eta = \eta^{-1} = \mbox{diag} (1,1,-1).
$$
Let $z=p+iq$, $w=r+is$, so $(p, q, r, s)$ are standard flat coordinates 
on $\R^{4}$. The gauge fields $A_{p},A_{q},A_{r},A_{s}$ are 
$\mathfrak{su}(2,1)$ valued. The relations $A_{z}=(A_{p}-iA_{q})/2$, 
$A_{\bar{z}}=(A_{p}+iA_{q})/2$ together with (\ref{liealgebra}) imply that
$$
\bar{A_{z}}^{t} = - \eta A_{\bar{z}} \eta^{-1},
$$
with a similar relation between $A_{w}$ and $A_{\bar{w}}$. 
Concretely, this means that
$$
A_{\bar{z}} = \left(
\begin{array}{ccc}
 a & b  & c  \\
 d & e  & f  \\
 g & h  & k  
\end{array}
\right), \ \ \ 
A_{z} = \left(
\begin{array}{ccc}
-\bar a  & -\bar d  & \bar g  \\
-\bar b  & -\bar e  & \bar h  \\
\bar c   & \bar f   & -\bar k  
\end{array}
\right),
$$
where $a+e+k=0$ (and of course $A_{w}$ and $A_{\bar{w}}$ are related 
in the same way).

Choosing a real form
$SU(2,1)$ of $SL(3,\C)$ on
restriction to the Euclidean slice imposes a constraint $\t U = \b
U$ and yields the affine sphere equation (\ref{LYZeq}). 
 
To sum up, one could achieve the characterisation of the ansatz (\ref{MDLYZ}), with
$\t z = \b z,$ $\t U= \b U,$
analogous to
Proposition \ref{gaugecond}.  Let us again choose the double 
null coordinates such
that the generators of the symmetry group of the 
ASDYM are given by $\p_{\t w}, \;
\p_w.$  With the chosen reality condition
the ASDYM equations reduce to the $SU(2, 1)$ Hitchin equations
\begin{eqnarray} 
\label{HitchinLYZ1}
D_zA_w &=& 0 \\
\label{HitchinLYZ2}
F_{z \bar z} + [A_w, A_{\bar w}] &=& 0,
\end{eqnarray}
where 
\be \label{bar} A_{\bar z} = - \eta^{-1} \bar{{A_z}}^{t} 
\eta \quad \mbox{and} \quad A_{\bar w} =
- \eta^{-1} \bar{{A_w}}^{t} \eta. \ee
We now consider the reduction of the system
(\ref{HitchinLYZ1}),(\ref{HitchinLYZ2}).
Theorem \ref{Main_thm_LYZ} arises as  a corollary of Proposition \ref{gaugecond}.\koniec

\subsection{Tziz\'eica equation}
The Tzitz\'eica equation (\ref{Teqn}) is a different real form of
(\ref{holTeqn}). It arises from the ASDYM 
with the gauge group 
$SL(3, \R)$ on restriction to the ultrahyperbolic real slice $\R^{2,2}$ in 
$\C^4$ with \[(w, \t w, x=z,
y=\t z )\] real.  The Higgs fields are given by $ P=A_{\t w}, Q=A_w$
and  the metric on the space of orbits of $X_1=\p_{\t w}$ 
and $X_2=\p_w$ has signature $(1, 1)$.

The real version of the 
ansatz (\ref{holMD_Higgs}),(\ref{holMD_potential}) can be characterised
analogously to the holomorphic case treated in Proposition \ref{gaugecond}. 
However, one
needs to take care of the fact that
$e^{u(x,y)} >0$ for real valued function $u(x,y).$  There are two places where
this needs to be considered.  First is where we use condition 
${\bf (i)}$ in
Proposition \ref{gaugecond} to put $(Q,P)$ in the form
(\ref{holMD_Higgs}),(\ref{holMD_potential}).  To write
$\mbox{Tr}(PQ)=e^{u(x,y)},$  we require that $\mbox{Tr}(PQ)>0$.
Assume that this can be done at a point $(x_0, y_0)$ (if not
then change coordinates $y\rightarrow -y)$ and restrict
the domain of $u$ to a neighbourhood of this point where
the positivity still holds.

The second place where the  problem of the sign arises is 
when we use the coordinate symmetry to transform
\[u_{xy} = e^u -r(x)b(y)e^{-2u}\]
to the Tzitz\'eica equation (\ref{Teqn}).  This can only be done for
$r(x)b(y)>0.$ 
The sign of $r(x)b(y)$ is governed by the quantity 
$ {\mbox {Tr}} \left( (D_xP)^2(D_{y}Q)^2 \right) $
in condition ${\bf (ii)}.$  To see it, note that in the notation of (\ref{gMDansatz}),
\[ {\mbox {Tr}} \left( (D_xP)^2(D_{y}Q)^2 \right) = (sktr) e^{2u}. \]
After we set $t=1,$ the condition
${\bf (iii)}$ implies that $k=e^u>0.$
Hence, the sign of $sr,$ and thus the sign $r(x)b(y)$ is the same as the
sign of $ {\mbox {Tr}} \left( (D_xP)^2(D_{y}Q)^2 \right).$  However, this
cannot be changed by real coordinate transformation $x \rightarrow \h
x(x), \; y \rightarrow \h y(y)$, 
because
\[ {\mbox {Tr}} \left( (D_xP)^2(D_{y}Q)^2 \right) \lra \left(
\frac{d{\h x}}{dx} \right)^{-4} \left( \frac{d{\h y}}{dy} \right)^{-4} 
{\mbox {Tr}} \left( (D_xP)^2(D_{y}Q)^2 \right), \]
where we have used $Q^2=0=P^2.$
Therefore, condition ${\bf (ii)}$ in Proposition \ref{gaugecond} needs to be replaced by
$ {\mbox {Tr}} \left( (D_zP)^2 \right) = 0 = {\mbox {Tr}} \left(
(D_{\t z}Q)^2 \right)$ and
\[{\mbox {Tr}}\left( (D_zP)^2(D_{\t z}Q)^2 \right) > 0 \] 
in the domain of $u$. 

\vspace{0.5cm}

\indent We remark that ${\mbox {Tr}} \left( (D_xP)^2(D_{y}Q)^2
\right) < 0$ corresponds to the equation
\[ u_{xy}= e^u + e^{-2u}, \]
whereas $ {\mbox {Tr}} \left( (D_xP)^2(D_{y}Q)^2 \right) =0$ 
yields Louiville equation
\[u_{xy}=e^u.\] 
Therefore, the sign of 
${\mbox {Tr}} \left( (D_xP)^2(D_{y}Q)^2 \right)$ corresponds to the
sign of $\epsilon$ in (\ref{definitesph}).
\section{$\Z_3$ two dimensional Toda chain}
\setcounter{equation}{0}
As a byproduct of the proof of Proposition \ref{gaugecond}, we
find that, dropping condition ${\bf (iii)}$ in this proposition,
the Hitchin equations can be
reduced to a coupled system which includes the $\Z_3$ two dimensional Toda chain
  \cite{Mikhailov} as a special case.  Recall that a two dimensional
  Toda chain is given by
\be \label{TodaN} (u_\a)_{xy} - e^{(u_{\a +1} - u_\a)} + e^{(u_{\a} -
  u_{\a-1})} =0,  \ee
where $\a \in \Z.$  In this paper (\ref{TodaN}) is called the $\Z_3$ two
dimensional Toda chain when \\
$ i) \quad \a \in \Z/ \Z_3 \qquad \mbox{and} $ \\
$ ii)\quad  u_1+u_2+u_3 =0.$ \\
We summarise the result in the following proposition.

\begin{prop} \label{Todaprop}
Let $u_1, u_2$ be functions of $(x,y).$
The coupled system of equations
\begin{eqnarray} \label{Toda}
(u_1)_{xy} - \epsilon_1 e^{(u_2-u_1)} + e^{2u_1+u_2} &=& 0 \nonumber \\
 (u_2)_{xy} + \epsilon_1 e^{(u_2-u_1)} - \epsilon_2 e^{-2u_2-u_1} &=& 0,
\end{eqnarray}
where $\epsilon_1, \epsilon_2 = \pm 1,$ is gauge equivalent to the 
$SL(3,\R)$  Hitchin equations {\em {(\ref{holHitchin1}, b, c)}}
with $z=x, \tilde{z}=y$ real, and \\
${\bf (i)}$ the Higgs fields $P$ and $Q$ have minimal polynomial $t^2,$ with ${\mbox
Tr}(PQ) \ne 0,$ \\
${\bf (ii)} \quad {\mbox {Tr}} \left((D_xP)^2 \right) = 0 = {\mbox Tr} \left(
(D_{y}Q)^2 \right)$ and ${\mbox {Tr}} \left( (D_xP)^2(D_{y}Q)^2
\right) \ne 0.$
\end{prop} 
{\bf Proof.}  These conditions are the first two conditions in
  Proposition \ref{gaugecond}.
 Following the proof  and assuming 
condition ${\bf(i)}$ gives (\ref{gaugecompat}).   
However, now it is not possible to use gauge symmetry
to set the diagonal elements of both $A_x$ and $A_y$ to be the same as
in (\ref{holMD_potential}) without the
compatibility condition.  Instead, let us use only the gauge transformation
(\ref{gdiag}) to eliminate the diagonal elements of $A_{y},$ by
choosing $(\ln a)_{y} = -p.$  

As before, condition ${\bf (ii)}$ implies that $m=h=0$ and $sktr \ne 0.$  The
Hitchin equations {\mbox {(\ref{holHitchin1}, b, c)}} imply that $t$ is a function
of $x$ only.  Hence, we can use the residual gauge freedom
(\ref{gdiag}) with $a=a(x)$ to set $t=1.$  Equation
(\ref{holHitchin3}) then gives
\begin{eqnarray} 
\label{thirdeqn1}
n_{y} + r(x)s &=& e^u  \\
\label{thirdeqn2}
2n_{y}-u_{xy}+r(x)s-k &=& 0 \\
\label{thirdeqn3}
 s_x + 3ns - su_x &=& 0 \\
\label{thirdeqn4}
 k_x +2ku_x - 3kn &=& 0. 
\end{eqnarray}
Equations (\ref{thirdeqn3}) and (\ref{thirdeqn4}) imply that $sk = c(y)e^{-u},$ where
$c(y)$ is some arbitrary function which arises from the integration.
Now, since $s \ne 0,$ let us write 
\[ k=\frac{c(y)}{s}e^{-u} \quad \mbox{and} \quad n= \a_x, \; s= \pm
e^\beta, \]
for some functions $\a(x,y)$ and $\beta(x,y).$
Then, (\ref{thirdeqn3}) becomes
\[e^\beta(\beta_x + 3 \a_x - u_x) =0,\]
which can be integrated to give
\[ s=b(y)e^{u -3 \a} \quad {\mbox{and}} \quad n= \a_x \]
for some $b=b(y) \ne 0.$
 Finally, (\ref{thirdeqn1}) and (\ref{thirdeqn2}) give a coupled system 
\begin{eqnarray} \label{coupled}
\a_{xy}+r(x)b(y)e^{u - 3\a} - e^u &=& 0  \nonumber \\
2\a_{xy} - u_{xy} +r(x)b(y)e^{u-3\a} - c(y)b^{-1}(y)e^{-2u +  3\a}
&=& 0.
\end{eqnarray}
Set $u_1= \a, \; u_2= -2 \a +u,$ and change the coordinate $y \rightarrow
-y.$ The system (\ref{coupled}) becomes
\begin{eqnarray*} 
(u_1)_{xy}-r(x)b(y)e^{u_2 - u_1} + e^{2u_1+u_2} &=& 0  \nonumber \\
(u_2)_{xy} + r(x)b(y)e^{u_2-u_1} - c(y)b^{-1}(y)e^{-2u_2 -u_1}
&=& 0,
\end{eqnarray*}
which can be transformed into (\ref{Toda}) by 
the change of dependent variables and coordinates.  There are four
distinct cases depending on the signs of
$\epsilon_1, \epsilon_2.$  Since the coordinates are real, the signs of
$\epsilon_1, \epsilon_2$ are the same as those of $r(x)b(y)$ and
$c(y)b^{-1}(y),$ respectively.   Similar to 
the real version of Proposition \ref{gaugecond} for the Tzitz\'eica
equation, $r(x)b(y)$ and
$c(y)b^{-1}(y)$ can be related to some gauge invariant quantities.  
It can be shown that at a given point $(x_0,y_0)$ the signs of $r(x)b(y)$ and
$c(y)b^{-1}(y)$ are determined by the signs of 
\begin{eqnarray*} 
{\bf (a)} &:=& {\mbox {Tr}}\left( (D_xP)^2(D_{y}Q)^2
\right), \\
{\bf (b)} &:=&
{\mbox {Tr}}\left( (PQ)^2(D_xP)(D_{y}Q) -
PQ(D_xP)QP(D_{y}Q) \right). 
\end{eqnarray*}
We shall analyse these signs and then restrict the domains of
$(u_1,u_2)$ to a neighbourhood of $(x_0,y_0)$ where the signs remain constant.
If ${\bf (a)}>0,$ setting $t=1$ gives $skr>0,$ which gives
 $r(x)c(y)>0.$  This implies that $r(x)b(y)$ and
$c(y)b^{-1}(y)$ have the same signs.  Now if ${\bf (b)}>0,$ then $k>0$
meaning $c(y)b^{-1}(y)>0,$ hence $r(x)b(y)>0.$  Similarly if ${\bf (b)}<0$
 then $c(y)b^{-1}(y)$ and $r(x)b(y)\; <0.$  On the other hand,
 ${\bf (a)}<0$ implies that $r(x)b(y)$ and
$c(y)b^{-1}(y)$ have opposite signs.  Then, the sign of ${\bf (b)}$
 determines the sign of $c(y)b^{-1}(y).$  The important point is that
 the signs of ${\bf (a)}$ and
 ${\bf (b)}$ cannot be changed by real coordinate transformations.
 This completes the proof.
\koniec

\section{Other gauges}
\setcounter{equation}{0}
There are several gauge inequivalent ways 
to reduce the ASDYM equations to
the  Tzitz\'eica equation or to the definite affine sphere equation. 
The reductions are relatively easy to obtain, but their gauge invariant characterisation
requires much more work. Here we shall mention
one other possibility which is not gauge equivalent to 
(\ref{holMD_Higgs}, \ref{holMD_potential}).

It can be shown that the holomorphic Tzitz\'eica equation
(\ref{holTeqn}) also arises
from the Hitchin equations with
\begin{eqnarray}
\label{Wang_fields}
P &=& \left(
\begin{array}{ccc}
0  & 0  & 1  \\
1 & 0  & 0  \\
0 & 1  & 0  
\end{array}
\right), \quad
Q = \left(
\begin{array}{ccc}
0  & e^{-2u}  & 0  \\
0 & 0  & e^u \\
e^u & 0  & 0  
\end{array}
\right), \nonumber \\
A_z &=& \left(
\begin{array}{ccc}
u_z  & 0  &  0 \\
0 & -u_z  & 0  \\
0 & 0  & 0  
\end{array}
\right), \quad
A_{\tilde{z}} = 0.
\end{eqnarray}
The real version of this ansatz was implicitly used by E. Wang \cite{Wang}.
  
 Let us 
comment on how this  formulation is  related to (\ref{holMD_Higgs}), (\ref{holMD_potential}).
First note that the Lax pairs (\ref{lax_pair_H}) with 
(\ref{holMD_Higgs}),(\ref{holMD_potential}) and (\ref{Wang_fields})
are equal for $\lambda=1$. Now consider the ansatz
(\ref{holMD_Higgs}),(\ref{holMD_potential}) and 
set $\lambda=1$ in the Lax pair
(\ref{lax_pair_H}). 
Introduce the 
new spectral parameter by exploiting the Lorentz symmetry and 
rescaling the coordinates 
\[
(z, \tilde{z})\longrightarrow (\hat{\lambda} z, \hat{\lambda}^{-1} \tilde{z})
\]
and read off
new $A_z, A_{\tilde{z}}, P, Q$ from  (\ref{lax_pair_H}) with $\lambda$ replaced by
$\hat{\lambda}$. This yields
the ansatz (\ref{Wang_fields}).

Choosing the Euclidean reality conditions and reducing the gauge group
to $SU(2, 1)$ we find another reduction of ASDYM to the affine sphere equation.
Take the following ansatz, in which the gauge 
fields are independent of $w$ and $\bar{w}$, $\psi=\psi(z, \bar{z})$ is a real 
function, and $U(z, \bar{z})$ is a complex function:
\begin{eqnarray}
\label{first_ansatz}
A_{w} &=& \left(
\begin{array}{ccc}
 0 & 0  & \frac{1}{\sqrt{2}} e^{\psi/2}  \\
 \bar{U} e^{-\psi} & 0  & 0  \\
 0 & \frac{1}{\sqrt{2}} e^{\psi/2}  & 0  
\end{array}
\right), \nonumber\\
A_{\bar{w}} &=& \left(
\begin{array}{ccc}
 0 & - U e^{-\psi}  & 0  \\
 0 & 0  & \frac{1}{\sqrt{2}}e^{\psi/2}  \\
 \frac{1}{\sqrt{2}} e^{\psi/2} & 0  & 0  
\end{array}
\right),\nonumber
\end{eqnarray}
\begin{eqnarray}
A_{z} &=& \left(
\begin{array}{ccc}
 -\frac{1}{2} \psi_{z} & 0  & 0  \\
 0 & \frac{1}{2} \psi_{z}  & 0  \\
 0 & 0  & 0  
\end{array}
\right),\nonumber \\
A_{\bar{z}}&=&\left(
\begin{array}{ccc}
 \frac{1}{2} \psi_{\bar{z}} & 0  & 0  \\
 0 & - \frac{1}{2} \psi_{\bar{z}}  & 0  \\
 0 & 0 & 0  
\end{array}
\right).
\end{eqnarray}
Recall that $A_w = Q$ and $A_{\b w} = -P.$  
The equation $F_{zw}= 0$ is satisfied provided that
\[
U_{\bar{z}}=0,
\]
i.e. $U$ must be holomorphic. The second ASDYM equation
$ F_{z \bar{z}} + F_{w \bar{w}} =0$ is satisfied 
if and only if (\ref{LYZeq}) holds.
\section{Semi--Flat Calabi--Yau metric}
\setcounter{equation}{0}
In this  section we consider the semi--flat Calabi--Yau metric
constructed by Loftin, Yau and Zaslow, and obtain the local expression
of the metric explicitly in term of solution of the definite affine sphere equation.

Let us first recall the Simon--Wang approach to
affine spheres \cite{SW}.   Consider the 
parametrisation of an elliptic affine sphere
\[(z, \bar z)
\mapsto f=(f^1(z, \bar z),f^2(z, \bar z),f^3(z, \bar z) ) \in \R^3.\]
 The structure equations\footnote{For the elliptic affine sphere with affine mean
curvature set to 1, the shape operator is $S=I.$ Now, with the affine
 metric (\ref{affineh}), the affine normal chosen to
point inward from the surface is given by minus the position vector
$-f,$ and the structure equations (\ref{Gauss}) and (\ref{Weingarten}) become
\begin{eqnarray*}  D_X f_* (Y) &=& f_*(\nabla_X Y) +
  h(X,Y)(-f)  \\
D_X (-f) &=& -f_*(X).
\end{eqnarray*}
Note that we have abused the notation so that $f$ also denotes the immersion.}
 defining the affine sphere can be written as a linear
first order system of PDEs in $f, f_z$ and $f_{\bar z}$
\begin{eqnarray}
\label{structureeqn}
\frac{\p}{\p z}  \left(
\begin{array}{c}
f  \\
f_z  \\
f_{\bar z}  
\end{array}
\right) = \left(
\begin{array}{ccc}
0 & 1 & 0 \\
0 & \psi_z & U e^{-\psi} \\
-\frac{1}{2} e^{\psi} & 0 & 0
\end{array} \right)
\left(
\begin{array}{c}
f  \\
f_z  \\
f_{\bar z}  
\end{array}
\right), \\ \nonumber
\frac{\p}{\p \bar z}  \left(
\begin{array}{c}
f  \\
f_z  \\
f_{\bar z}  
\end{array}
\right) = \left(
\begin{array}{ccc}
0 & 0 & 1 \\
-\frac{1}{2} e^{\psi} & 0 & 0 \\
0 & \bar U e^{-\psi} & \psi_{\bar z}
\end{array} \right)
\left(
\begin{array}{c}
f  \\
f_z  \\
f_{\bar z}  
\end{array}
\right),
\end{eqnarray}
where we have set
the affine mean curvature to 1.
The compatibility condition for this over-determined system is the affine sphere
equation (\ref{LYZeq}).

Therefore, given a solution $\psi,$ one can find $f$
and hence the cone over the sphere
\be \label{f}
(z, \b z, r) \longmapsto (x^1=rf^1(z, \bar z), \; x^2=rf^2(z, \bar z),
\; x^3=rf^3(z, \bar z)).
 \ee
This expression can be inverted locally to give $r=r(x)$.
\vskip5pt

\noindent {\bf Proof of Proposition \ref{CYprop}.}
The metric cone over an elliptic affine sphere is given by  (\ref{gB})
with $\phi(x) = r^2/2$ and the corresponding semi-flat metric
(\ref{g_omega}).

The matrix $\phi_{jk}$ in (\ref{g_omega}) can be obtained by contracting the metric (\ref{gB}) with
$\p/\p x^j,\p/\p x^k.$
Given a solution of the affine sphere equation $\psi,$ we know $g_B$ in the basis
$(dr, dz, d \bar z),$
thus we want to express $\p/\p x^j$ in terms of
$\p/\p r, \p/\p z, \p/\p \bar z.$
Now, from (\ref{f}), we have that 
\[\left( \begin{array}{c}
\p/\p x^1  \\
\p/\p x^2  \\
\p/\p x^3  
\end{array} \right) = N^{-1}
 \left(
\begin{array}{c}
\p/\p r  \\
r^{-1} \p/\p z  \\
r^{-1} \p/\p \bar z  
\end{array}
\right), \quad
\mbox{where} \quad
 N=
\left(
\begin{array}{ccc}
f^1 & f^2 & f^3 \\
f^1_z & f^2_z & f^3_z \\
f^1_{\bar z} & f^2_{\bar z} & f^3_{\bar z}
\end{array} \right).
\]
Moreover, $N$ is the matrix solution of the linear
system (\ref{structureeqn}), whose existence and the existence of its inverse $N^{-1}$
are guaranteed by the affine sphere equation.
Writing 
\[ N^{-1}=
\left(
\begin{array}{ccc}
p_1 & q_1 & \b q_1 \\
p_2 & q_2 & \b q_2 \\
p_3 & q_3 & \b q_3
\end{array} \right),
\]
one calculates $\phi_{jk}$ and thus the metric on the fibre to be
\[\phi_{jk}dy^jdy^k = (p_j p_k + e^\psi q_j \b q_k)dy^jdy^k.\]
Now, let us introduce new coordinates
\[ \tau:=p_iy^i, \quad  \xi:=q_iy^i, \quad \b \xi := \b q_iy^i  \]
and write $p_idy^i = d\tau-y^idp_i$ etc.  Denote the two matrices of
coefficients in the linear system (\ref{structureeqn}) by $-A^{(z)}$ and $-A^{(\b z)}$
respectively, so that (\ref{structureeqn}) is
\[ \p_z N + A^{(z)} N=0, \qquad  \p_{\b z} N + A^{(\b z)} N=0. \]
Then, by considering the corresponding equation for $N^{-1},$
the one--forms $y^idp_i, y^idq_i, y^id\b q_i$ can be written
in terms of coordinates $\tau, \xi, \b \xi$ and components of $A^{(z)}$ and $A^{(\b z)},$
which are known in terms of $\psi$.  

Finally, we can write the
metric (\ref{g_omega}) as
\[ g = dr^2 + r^2e^\psi |dz|^2 + |d\tau+ \a|^2 + e^\psi |d\xi+ \beta|^2,
\]
where 
\[
\a = -\frac{1}{2}e^\psi(\b \xi dz + \xi d\b z), \quad
\beta = (\tau+\xi \psi_z)dz + e^{-\psi}\b U \b \xi d\b z.
\]
By similar calculation, the K\"ahler form can be written as
\[ \om = dr \wedge (d\tau+\a) + \frac{r}{2}e^{\psi}(d\b z \wedge
(d\xi +\beta) +dz \wedge (d\b \xi+ \b \beta)  ).\]

Using the relation between the metric, the K\"ahler form and the
complex structure, we find holomorphic basis $\{e_1, e_2, e_3\}$ (\ref{basis}) and
write $g$ and $\om$ as in Proposition \ref{CYprop}, where we have
introduced a complex coordinate $w=r+i\tau.$ 
\koniec
{\bf Remark 1.}  The Ricci flat condition for
the metric (\ref{semiflatcymetric}) reduces to the affine sphere
equation (\ref{LYZeq}) for
$\psi(z, \b z)$ and $U(z)$. Equation (\ref{LYZeq}) is invariant under 
the transformations 
$ \p/\p z \ra  \p/\p {\h z},\; \psi \ra \h \psi,\; U \ra \h U,$ where
\[ \p/\p_{\h z} = e^{-j(z)}\p/\p_z, \quad  \h \psi =
\psi -j(z) -\overline{j(z)}, \quad \mbox{and} \quad \h U = e^{-3j(z)}U.
\]
This can be understood geometrically, as $e^\psi dz d\b z$ and $U
dz^3$ are the affine metric and the cubic differential respectively of the
affine sphere.  The metric  (\ref{semiflatcymetric}) is invariant
under the above transformations, together with $\xi \ra \h \xi =
e^{j(z)} \xi.$

\vspace{0.5cm}

\noindent {\bf Remark 2.}  One expects the linear system associated 
with the structure
equations of affine spheres (\ref{structureeqn}) to be equivalent to
the Hitchin Lax pair (\ref{lax_pair_H}) giving rise to the affine
sphere equation.  The matrices
$A^{(z)}$ and $A^{(\b z)}$ in (\ref{structureeqn}) are unique up to 
gauge transformations
\[ A^{(z)} \lra g^{-1}A^{(z)}g + g^{-1}\p_z g, \quad  A^{(\b z)} \lra
g^{-1}A^{(\b z)}g + g^{-1}\p_{\b z} g.\]
If we write 
\be \label{A(z)} A^{(z)} =  (A_z + \lambda P), \quad A^{(\b z)} =  (A_{\b z} +
\lambda^{-1} Q) \ee
for some value of $\lambda,$ then it follows that $(A_z, A_{\b z}, Q, P)$ 
will satisfy the Hitchin equations
(\ref{holHitchin1}, b, c), with reality condition $\t z=\b z$.
Conversely, given a solution $(A_z, A_{\b z}, Q, P)$ to the Hitchin
equations, we should be able to find a value of spectral parameter
$\lambda$ such that $(A_z + \lambda P)$ and $(A_{\b z} +
\lambda^{-1} Q)$ can be gauge transformed to  $A^{(z)}$ and $A^{(\b
  z)}$ respectively.

For example, we can obtain $A^{(z)}$ and $A^{(\b z)}$ in
(\ref{structureeqn}) from the ansatz (\ref{MDLYZ}), with $\t z = \b z$
and $\t U = \b U,$ by gauge
transformation with 
\[ g =  \left( 
  \begin{array}{ccc}
    1 & 0 & 0 \\
    0 &-\sqrt{2} e^{-\psi/2} & 0 \\
    0 & 0 & -\sqrt{2} e^{-\psi/2}
  \end{array}
  \right),
 \]
and choosing the value of spectral parameter in (\ref{A(z)}) to be
$\lambda =1.$  Note that we need $\det g \ne 1,$ since $A^{(z)}$ and
$A^{(\b z)}$ are not traceless.

\section{Painlev\'e III}
\setcounter{equation}{0}
One of the main results of Loftin, Yau and Zaslow \cite{LYZ} is the
existence of radially symmetric solutions of the affine sphere equation
(\ref{LYZeq}) for $U(z)=z^{-2},$ with prescribed behaviour near the
singularity $z=0.$  In this section we shall show that the radially symmetric
solutions of (\ref{LYZeq}) are Painlev\'e III transcendents.

{\bf Proof of Proposition \ref{painleve_prop}.} Set $U=z^{-2}$, 
and look for solutions of (\ref{LYZeq}) of the 
form $\psi=\psi(\rho)$, where $\rho=|z|$. 
Making a substitution $\psi(\rho)=\log{(\rho^{-3/2}H(\rho))}$
and introducing a new independent 
variable by $\rho=s^2$ yields the the following 
ODE for $H=H(s)$
\be
\label{PIII}
H_{ss}=\frac{{(H_s)}^2}{H}-\frac{H_s}{s}-\frac{8H^2}{s} -\frac{16}{H}.
\ee
This is the celebrated Painlev\'e III equation \cite{Ince}
\[H_{ss}= \frac{(H_s)^2}{ H}-\frac{H_s}{s}+
\frac{\alpha H^2 +  \beta  }{s}+\gamma H^3 + \frac{\delta }{H}
\]
with special values of parameters
\[
(\alpha, \beta, \gamma, \delta)=(-8, 0, 0, -16).
\] 
In the classification of Okamoto \cite{Okamoto} it falls in
the type  $D7$.  
\koniec
{\bf Remarks.}
\begin{itemize}
\item
One can consider the radial symmetry reduction of the
 affine sphere equation (\ref{LYZeq}) with $U= z^{-n}$ for general $n
 \in \Z.$  
\begin{enumerate}
\item[${\bf n \ne 3.}$]  Changing the independent variable to
 \[s=(z \b z)^\frac{3-n}{4}\] and using the ansatz \[\psi = \log
 \left(s^{-\left(\frac{1+n}{3-n}\right)} H(s)^k \right)\] 
with $k=\pm 1$
reduces 
(\ref{LYZeq}) to the Painlev\'e III equation with
 parameters $(\a, \beta, \gamma, \delta) = \left( \frac{-8}{(3-n)^2},
 0, 0, \frac{-16}{(3-n)^2} \right)$ and 
$(\a, \beta, \gamma, \delta) = \left( 0, \frac{8}{(3-n)^2},
 \frac{16}{(3-n)^2}, 0 \right)$ for $k=1$ and $k= -1,$ respectively.  In both cases, the
 Painlev\'e III equations are of type $D7$ in Okamoto's classification.  
\item[{\bf n=3.}]
Setting  $\psi=\psi(s)$ where
$s=(z \b z)^\frac{1}{2}$ in equation (\ref{LYZeq}) yields
\be \label{radialn3} \psi_{ss}+ \frac{\psi_s}{s} +\frac{4 e^{-2\psi}}{s^6}+2e^\psi =
 0,\ee 
which, under multiplication by
$\left( \frac{\psi_s}{2}+\frac{1}{s} \right),$ gives a first-order ODE
\be \label{1stODE} \frac{{\psi_s}^2}{4}+ \frac{\psi_s}{s} + e^\psi
 -\frac{e^{-2\psi}}{s^6} +\frac{c}{s^2} = 0,   \ee
where $c$ is a constant of integration.  Hence any solution to (\ref{1stODE}) 
such that $s\psi_s\neq -1/2$ gives rise to  a solution to
(\ref{radialn3}), and conversely all solutions to (\ref{radialn3})
arise from (\ref{1stODE}). 
Equation (\ref{1stODE})  is integrable by quadratures in terms of
the elliptic functions. 
\end{enumerate}
\item
In general, a Painlev\'e III equation may have two types
of special (i.e. non--transcendental) solutions: the finite
number of rational solutions
and a one parameter family of Riccati type solutions expressible by
special functions \cite{Ince}. For the values of parameters in (\ref{PIII})
the Riccati solutions do not exist, and there exists a unique
algebraic solution
\[
H=-(2s)^{1/3}.
\]
This corresponds to 
\[
\psi=\frac{1}{3}\log{(2)}-\frac{4}{3}\log{(|z|)}+\log{(-1)}
\]
which is not real.
There are B\"acklund transformations leading to new solutions, but
they change the value of the parameters. This shows that the desired 
radial solution to the affine sphere equation (\ref{LYZeq}) is
transcendental.
In \cite{K, BE} it has been shown 
that the radial solutions of the Tzitz\'eica equation (\ref{Teqn}) also
satisfies Painlev\'e III of type $D7.$
\end{itemize}

\subsection{Lax pair for Painlev\' e III}
The standard isomonodromic approach to Painlev\'e III identifies this equation
with $SL(2, \C)$ isomonodromic problem with two double
poles. The connection with affine differential geometry and its underlying 
isospectral Lax pair  suggests that there is an alternative
isomonodromic Lax pair for PIII given in terms of 3 by 3 matrices,
as opposed to the standard Lax pair
with 2 by 2 matrices \cite{JM}. 
(See also \cite{MW} where $SL(2, \C)$ ASDYM has been reduced to PIII).

Let us now return to the holomorphic setting, and consider the Lax
pair for ASDYM in $\C^4$ with gauge group $SL(3,\C)$
\[  (D_w + \ll D_{\t z})\Psi =0, \quad (D_z + \ll D_{\t w})\Psi = 0, \]
where $\Psi$ is a vector-valued function of $w, \t w, z, \t z$ and $\lambda.$
We require that the connection is invariant under the 3 dimensional
subgroup of the conformal group $PGL(4,\C)$ generated by 
\be \label{generators} 
\{ \p_w, \; \p_{\t w}, \; z\p_z - {\t z}\p_{\t z} \}, \ee
and introduce coordinates $(\rho, \theta) \in \C^2$ such that $z=\rho e^{i
  \theta},$ $\t z = \rho e^{-i \theta},$ and $z\p_z - {\t z}\p_{\t z}=
  -i\frac{\p}{\p \theta}.$  Then the ASDYM Lax pair becomes
\begin{eqnarray*}
\left( -\zeta \p_\rho + \rho^{-1}\zeta^2\p_\zeta + 2(A_w - \zeta
  e^{-i\theta}A_{\t z}) \right)\Psi&=&0, \\
\left( \p_\rho + \rho^{-1}\zeta \p_\zeta + 2 (e^{i\theta}A_z - \zeta A_{\t w})
  \right) \Psi &=&0, 
 \end{eqnarray*}
where  the
gauge fields are in an invariant gauge; $(A_w, A_{\t w},
e^{i\theta}A_z, e^{-i\theta}A_{\t z})$ are functions of $\rho$ only, 
and $\zeta= -\ll e^{i\theta}$ is an invariant spectral
parameter\footnote{The spectral parameter $\ll$ is not constant along
  the lift of the generators (\ref{generators}) to $\C^4 \times \CP^1
  \in (w, \t w, z, \t z, \ll)$ where $\Psi$ is defined. However, the
  invariant spectral parameter $\zeta$ is constant along the lift, and
  hence we are allowed to express $\Psi$ as a function of $\rho$ and
  $\zeta$ only.}.
Taking linear combinations of these two linear PDEs gives a Lax pair of the form
\be \label{reducedlax} \frac{\p \Psi}{\p \zeta} = \hat L \; \Psi, \quad \frac{\p
  \Psi}{\p \rho}= \hat M \; \Psi,\ee
where 
\begin{eqnarray*}
\hat L &=& \rho \zeta^{-2} \left(\zeta^2A_{\t w}-A_w +\zeta(e^{-i\theta}A_{\t z} -
e^{i\theta}A_z) \right) \\
\hat M &=& \zeta^{-1}\left(A_w + \zeta^2A_{\t w} - \zeta(e^{i\theta}A_z +
e^{-i\theta}A_{\t z}) \right).
\end{eqnarray*}
The calculation leading to Painlev\'e III (\ref{PIII}) implies that
if we gauge transform ansatz (\ref{MDLYZ}) with
  $U(z)=z^{-2},$ $\t U(\t z) = \t z^{-2}$ into an invariant gauge and
substitute it into (\ref{reducedlax}), then in the new coordinate $s=\rho^{1/2}$ the
system (\ref{reducedlax}) becomes Lax pair of the Painlev\'e III with special
values of parameters (\ref{PIII}).  We shall now present this calculation:

An invariant gauge of (\ref{MDLYZ}) can be obtained using the gauge
transformation with
\[ g =  \left( 
  \begin{array}{ccc}
    e^{i \theta/3} & 0 & 0 \\
    0 & e^{-i2\theta/3} & 0 \\
    0 & 0 & e^{i\theta/3}
  \end{array}
  \right),
 \]
which does not change $A_w$ and $A_{\t
  w},$ but gives
 \begin{eqnarray*} e^{i\theta}A_{z} &=& \left( 
  \begin{array}{ccc}
    \frac{1}{6\rho} & \frac{1}{\sqrt{2}}e^{\frac{\psi}{2}} & 0 \\
    0 & -\left( \frac{1}{4}\psi_\rho + \frac{1}{3\rho} \right) & -\frac{1}{\rho^2}e^{-\psi} \\
    0 & 0 & \frac{1}{4}\psi_\rho + \frac{1}{6\rho}
  \end{array}
  \right), \\
 e^{-i\theta}A_{\t z} &=& \left( 
  \begin{array}{ccc}
    -\frac{1}{6\rho} & 0 & 0 \\
    -\frac{1}{\sqrt{2}}e^{\frac{\psi}{2}} & \frac{1}{4}\psi_\rho + \frac{1}{3\rho}  & 0 \\
    0 & -\frac{1}{\rho^2}e^{-\psi} & -\left(\frac{1}{4}\psi_\rho +
    \frac{1}{6\rho} \right)
  \end{array}
  \right).
\end{eqnarray*} 
Then, in terms of $s = \rho^{1/2}$ and $H(s)= s^3 e^\psi,$ the system
(\ref{reducedlax}) gives a Lax
pair for the Painlev\'e III equation (\ref{PIII}) as
\be \label{PIIIlax} \frac{\p \Psi}{\p \zeta} =  L \; \Psi, \quad \frac{\p
  \Psi}{\p s}=  M \; \Psi,\ee
where 
\begin{eqnarray*} L &=& -\frac{1}{\zeta^2}\left( 
  \begin{array}{ccc}
    \frac{\zeta}{3} & \frac{1}{\sqrt{2}}\zeta (sH)^{1/2} & \frac{1}{\sqrt{2}}(sH)^{1/2} \\
     \frac{1}{\sqrt{2}}\zeta (sH)^{1/2} & \zeta \left(\frac{1}{12} -\frac{sH_s}{4H}
      \right) & -\zeta \frac{s}{H} \\
     \frac{1}{\sqrt{2}}\zeta^2 (sH)^{1/2} & \zeta \frac{s}{H} &
     \zeta \left(\frac{sH_s}{4H}- \frac{5}{12} 
     \right)
  \end{array}
  \right), \\
M &=& \sqrt{2} \left( \frac{H}{s} \right)^{1/2} \left( 
  \begin{array}{ccc}
    0 & -1 & \frac{1}{\zeta} \\
    1 & 0  & \sqrt{2} \left(\frac{s}{H^3}\right)^{1/2}\\
    -\zeta & \sqrt{2} \left(\frac{s}{H^3}\right)^{1/2} & 0 
  \end{array}
  \right).
\end{eqnarray*} \\
The matrix $L$ has two double poles as expected for Painlev\'e III \cite{JM}, 
at $\zeta = 0$ and $\zeta =
\infty.$  

We note here that a different (i.e. gauge inequivalent) $3 \times 3$ isomonodromic Lax pair for Painlev\'e
III of type $D7$ was used
by Kitaev in \cite{K}.  The Lax pair can also be derived from the ASDYM
Lax pair, from a solution to Hitchin
equations which is gauge equivalent to (\ref{Wang_fields}).

\section*{Acknowledgements}
We  wish to thank Philip Boalch, Robert Conte, Eugene Ferapontov, Nigel Hitchin, John Loftin, Ian
McIntosh, Yousuke Ohyama and Wolfgang Schief for valuable comments.  
Prim Plansangkate is grateful to the
Royal Thai Government for funding her research.

\end{document}